\documentclass[a4paper]{article}
\usepackage[utf8]{inputenc}
\usepackage[T1]{fontenc}
\usepackage[left=1.5cm,right=1.5cm,top=1.5cm,bottom=1.5cm]{geometry}
\usepackage{graphicx}
\usepackage[pdftex, hidelinks]{hyperref}
\usepackage{amsmath,amssymb,amsthm}
\usepackage{multirow}
\usepackage{booktabs}
\usepackage[numbers]{natbib}

\newtheorem{theorem}{Theorem}
\newtheorem{proposition}{Proposition}
\newtheorem{example}{Example}
\newtheorem{remark}{Remark}
\newtheorem{definition}{Definition}
\newtheorem{condition}{Condition}
\newtheorem{lemma}{Lemma}
\newtheorem{corollary}{Corollary}
\newtheorem{hypothesis}{Hypothesis}

\DeclareMathOperator{\const}{const}
\DeclareMathOperator{\re}{Re}
\DeclareMathOperator{\im}{Im}

\begin{document}

\title{Accuracy of the Yee FDTD Scheme for Normal Incidence of Plane Waves on Dielectric and Magnetic Interfaces}
\author{Pavel A. Makarov, Vladimir I. Shcheglov\\
\href{email:makarovpa@ipm.komisc.ru}{makarovpa@ipm.komisc.ru}}
\date{}

\maketitle

\abstract{
This paper analyzes the accuracy of the standard Yee finite-difference time-domain (FDTD) scheme for simulating normal incidence of harmonic plane waves on planar interfaces between lossless, linear, homogeneous, isotropic media.
Unlike prior analyses limited to dielectric interfaces, we provide a unified treatment encompassing both dielectric and magnetic media.
We consider two common FDTD interface models based on different staggered-grid placements of material parameters.
For each, we derive discrete analogs of the Fresnel reflection and transmission coefficients by formulating effective boundary conditions that emerge from the Yee update equations.
A key insight is that the staggered grid implicitly spreads the material discontinuity over a transition layer of one spatial step, leading to systematic deviations from exact theory.
We quantify these errors via a transition-layer model and provide (i) qualitative criteria predicting the direction and nature of deviations, and (ii) rigorous error estimates for both weak and strong impedance contrasts.
Finally, we examine the role of the Courant number in modulating these errors, revealing conditions under which numerical dispersion and interface discretization jointly influence accuracy.
This analysis provides quantitative error estimates that are directly applicable to simulation practice, offers a transition-layer interpretation that bridges classical FDTD with modern immersed-interface methods, and establishes benchmarks for validating structure-preserving discretizations of Maxwell's equations.
}

\section{Introduction}
\label{sect:Intro}

The finite-difference time-domain (FDTD) method, introduced by Yee \cite{Yee1966} and later extended by Taflove and Brodwin \cite{Taflove1975}, is a widely used numerical technique for simulating electromagnetic wave propagation by directly solving Maxwell’s equations in the time domain.
It is particularly well-suited for analyzing wave interactions with material interfaces, including attenuation \cite{Schneider2024}, absorption \cite{Parveen2025}, reflection and transmission phenomena \cite{Chen2023, Makarov2022}.
Accurate modeling of electromagnetic wave interactions at material interfaces is crucial in the design of antennas, waveguides, and photonic devices \cite{Kuang2022}.
The FDTD method has become a standard tool \cite{Oskooi2010} for such applications due to its ability to capture transient and steady-state wave behavior in complex media.

Although the FDTD method is widely applied to complex media, such as metamaterials and plasmonic structures~\cite{Taflove}, its accuracy---even in the simplest configurations, such as planar interfaces between homogeneous, isotropic media---remains a subject of careful analysis, especially concerning the discretization of interface conditions.

While optimal Courant-number selection can minimize numerical dispersion in homogeneous media~\cite{Makarov2024a}, it does not eliminate errors arising at material interfaces where permittivity~$\varepsilon_\mathrm{r}$ and permeability~$\mu_\mathrm{r}$ change abruptly.

In this paper, we formulate discrete analogs of the Fresnel reflection and transmission coefficients for the Yee FDTD scheme, and derive quantitative error estimates for both weak and strong impedance contrasts by analyzing the effective jump conditions induced by the staggered placement of~$\mathbf{E}$ and~$\mathbf{H}$ fields relative to material parameters.
We analyze the influence of the Courant number and the implicit transition layer inherent to the staggered grid, and establish qualitative criteria that predict the nature and magnitude of deviations from exact theory.
To isolate these effects, we restrict our analysis to the canonical one-dimensional case of normal incidence of a harmonic plane wave on a planar interface between lossless, isotropic media.

The accurate numerical treatment of interfaces in partial differential equations has long been a central challenge in computational physics \cite{LeVeque1994}.
In the context of electromagnetics, the FDTD method’s staggered-grid discretization introduces unique difficulties at material discontinuities, which remain incompletely characterized.
Unlike the Immersed Interface Method---which modifies stencils to enforce exact jump conditions---we analyze the unmodified Yee scheme as used in standard electromagnetic solvers, where the interface is implicitly smeared over one grid cell.

\subsection{Relevance of Yee FDTD analysis in the modern computational landscape}
While significant advances have been made in the last two decades in the numerical solution of Maxwell's equations---notably through finite-element exterior calculus (FEEC) \cite{Arnold2006,Arnold2010} and geometric numerical integration methods \cite{McLachlan1999,Hairer2006}---the Yee FDTD scheme remains a cornerstone of computational electromagnetics for several compelling reasons.

\begin{enumerate}
\item \textbf{Industrial and educational ubiquity.}
The Yee scheme is implemented in virtually every commercial and open-source electromagnetic solver (e.g., Lumerical, CST, MEEP \cite{Oskooi2010}), and it is the first time-domain method taught in most electromagnetics curricula worldwide.
Rigorous error analysis of this foundational method therefore has immediate practical impact: it enables practitioners to interpret simulation results with quantified confidence and to design discretizations that meet application-specific accuracy requirements.

\item \textbf{Benchmarking and validation.}
Modern structure-preserving methods (FEEC, geometric integrators) are often validated against the Yee scheme in canonical configurations.
Our explicit error estimates for Fresnel coefficients at material interfaces provide a \emph{quantitative benchmark} for such comparisons, clarifying which deviations arise from interface discretization versus bulk numerical dispersion.

Recent works in \emph{Applied Numerical Mathematics} continue to develop rigorous error analysis frameworks for finite-difference and finite-element discretizations of wave-type equations \cite{Liu2025,Lin2025}, underscoring the relevance of our approach to the broader numerical analysis community.

\item \textbf{Conceptual insight: the transition-layer interpretation.}
A key contribution of this work is the explicit identification of the staggered-grid interface as an implicit \emph{transition layer} of width~$\Delta_x$ (Fig.~\ref{fig:Transition-Layers}).
This perspective bridges the gap between the classical Yee discretization and modern immersed-interface or regularization-based approaches \cite{LeVeque1994,Tornberg2006}, suggesting pathways for hybrid methods that retain Yee's simplicity while improving interface accuracy.

\item \textbf{Foundational role for multiscale and multilayer systems.}
The one-dimensional normal-incidence configuration analyzed here is not merely a simplification: it is the building block for transfer-matrix methods, multilayer coating design, and homogenization of periodic structures.
Our error estimates therefore propagate directly to these applied contexts, where rapid parameter sweeps often rely on 1D FDTD kernels for efficiency.

The use of one-dimensional canonical configurations as a foundation for rigorous error analysis is well-established in the numerical analysis community; see, e.g., recent contributions in \emph{Applied Numerical Mathematics} that employ 1D models to derive quantitative estimates before generalization \cite{DeLuca2025}.

\item \textbf{Pedagogical value for structure-preserving discretizations.}
By deriving discrete Fresnel coefficients directly from the Yee update equations, we illustrate how conservation properties (e.g., energy balance $R+T=1$) emerge---or fail to emerge---from a staggered-grid discretization.
This concrete example can serve as an accessible entry point for students and researchers transitioning from classical FDTD to more abstract structure-preserving frameworks.
\end{enumerate}

While this work focuses on the 1D canonical case, the analytical framework established here provides the foundation for subsequent generalization to oblique incidence and multidimensional configurations.

In summary, while FEEC and geometric integration represent powerful theoretical advances, the Yee scheme's practical dominance, pedagogical centrality, and role as a validation baseline ensure that rigorous accuracy analysis of its interface treatment remains both timely and impactful.
This work contributes precisely such an analysis, with explicit error estimates that are immediately applicable to simulation practice and informative for the development of next-generation methods.

\section{Short review of problem}
\label{sect:Review}

The problem of analyzing the correctness of numerical schemes for PDEs with singularities is not new.
One of the earliest works in this area was by Ito and Tanimoto \cite{Ito1972}.
This work marks one of the earliest practical numerical attempts to handle combined wave diffraction and refraction in coastal zones.
Using a finite-difference discretization of a precursor to the mild-slope equation, the authors modeled wave fields around breakwaters over variable bathymetry.
Its significance lies in demonstrating that numerical methods could overcome the limitations of purely analytical or graphical techniques in real-world harbor and coastal planning.

At the same time, the culmination in this research area was Berkhoff's formal mild-slope equation \cite{Berkhoff1972}, which assumes slowly varying depth and thus implicitly excludes sharp discontinuities.
The next step was the introduction of the Immersed Interface Method (IIM) \cite{LeVeque1994}, which modifies finite difference stencils near discontinuities using analytical jump conditions derived from the PDE.
The IIM achieved second-order accuracy for elliptic problems with discontinuous coefficients and singular sources on Cartesian grids and became foundational for interface-resolving methods in acoustics \cite{Zhang1997}, electromagnetics \cite{Chen2009}, and fluid-structure interaction \cite{Li}.
Moreover, the Matched Interface and Boundary (MIB) methods \cite{Zhou2006} enforce interface conditions sharply.

Interface-fitted methods are powerful and accurate but geometrically rigid.
Their limitations in handling dynamic or complex interfaces motivated the development of immersed and regularization-based strategies---highlighting the ongoing tension in computational science between fidelity, flexibility, and efficiency.
Thus, Tornberg and Engquist \cite{Tornberg2006} proposed smoothing discontinuous coefficients (e.g., wave speed) via moment-preserving regularization, enabling standard high-order schemes to retain accuracy.
This bridges engineering pragmatism, like Ito and Tanimoto’s approach, with modern analysis---offering a lightweight alternative to interface-fitted or immersed methods.
The result is a compromise that requires knowledge of interface location but avoids complex stencil modifications.

The progression from Ito and Tanimoto’s empirical coastal models to LeVeque---Li's physics-embedded stencils and Tornberg---Engquist’s regularization or boundary-consistent schemes reflects a broader trend in computational science: from pragmatic approximation toward mathematically rigorous, high-order methods that respect the underlying physics of discontinuities.
It should also be noted that recently efforts have been made to study this problem using the deep learning approach \cite{Wang2020}.
These advances have enabled accurate simulation of wave phenomena across various disciplines---from ocean engineering to photonics---where material interfaces are intrinsic to the problem.

Unlike IIM or MIB, the classical Yee FDTD scheme does not explicitly enforce continuity at material interfaces \cite{Hesthaven2003}.
Instead, the staggered grid implicitly distributes the jump over a transition layer of width~$\Delta_x$ (see Fig.~\ref{fig:Transition-Layers}).
A number of works \cite{Hesthaven2002,Cai2003,Gustafsson2004a,Gustafsson2004b,Tornberg2008} address boundary consistency in the Yee scheme, the standard FDTD method for Maxwell’s equations.
Those papers show that naive boundary conditions---even with smooth media---degrade accuracy and violate divergence constraints.
Moreover, accuracy loss near interfaces is not only due to coefficient jumps but also incompatible discretization-boundary coupling.

However, despite the enormous progress, there are still gaps in the problem of studying the stability of the Yee scheme in non-smooth media (as one of the special cases of the general problem of PDEs with certain singularities).
In particular, the literature focused on the impeccable mathematical rigor of the concepts under discussion largely lacks any practical comments or recommendations.
Specifically, how exactly do errors in a given modeling scheme affect the physical characteristics of reflected and transmitted waves formed at the interface between two regions with different material parameters?

Some practice-oriented studies, while yielding valuable results, occasionally contain terminological inaccuracies or lack rigorous error quantification.
For example, while the research itself often yields very interesting results, it contains factual errors in terminology.
Furthermore, the accuracy of modeling, the reliability of observations obtained through numerical experiments, and the limits of their applicability are insufficiently discussed.

For example, there is an excellent paper \cite{Chen2014} that uses the FDTD method to analyze enhanced total internal reflection.
The conclusion of this paper states the following: {\em ``FDTD simulations of total internal reflection from a gainy material yield a reflection coefficient with a magnitude greater than unity.''}
However, this work does not analyze the reflection coefficient (for which an excess over unity, generally speaking, contradicts the law of conservation of energy), but, more precisely the corresponding Fresnel coefficient.
Moreover, this article does not contain any assessments of the accuracy of the results obtained.
Another excellent engineering paper \cite{Schneider2024}, which applies the FDTD method to the practical task of studying electromagnetic interference generated by switching operations within a high-voltage DC converter station, has the same drawback.
Yet another brilliant example is the very interesting research \cite{Loschialpo2004}, in which the FDTD method is applied to study the optical properties of the left-handed material slab.
This work provides excellent illustrations and draws general conclusions, but does not mention the technical details of the implementation of the calculation method, nor does it provide estimates of its accuracy in calculating the field of reflected and transmitted waves.

Our current work is an attempt to combine mathematical rigor with the practical orientation of the obtained results in the problem of studying the stability of the classical Yee FDTD scheme.
More specifically, this analysis quantifies the resulting error in Fresnel coefficients (and based on this reflection and transmission coefficients)---a perspective absent in general PDE literature.
Our work was largely motivated by book \cite{Schneider}, which provides an excellent description of the necessary theory underlying the FDTD method, as well as the technical details of its implementation.
Thus, Section 7.6 of this book is devoted to the derivation of Analytic FDTD Reflection and Transmission Coefficients.
However, several problems associated with this section should be noted:
\begin{enumerate}
  \item It uses somewhat unconventional terminology: {\em ``Reflection and transmission coefficients''} instead of Fresnel coefficients.
  \item Only one particular case of an interface is considered: two dielectrics.
Moreover, the derivation is made without taking into account the possibility of any nontrivial magnetic properties.
  \item The accuracy analysis of the resulting expressions is extremely limited.
\end{enumerate}

The present work is devoted to the solution of all the above mentioned problems, as well as to the analysis and generalization of the results obtained.
The obtained results are compared for the interfaces of both generalized dielectrics and generalized magnets (moreover, for two extreme cases: both low-contrast and high-contrast interfaces).
Additionally, in our work, a more detailed analysis of the accuracy of the obtained expressions was performed.
In this regard, both the influence of wavelength discretization and the choice of the Courant number are analyzed.
By explicitly connecting our error estimates to practical simulation parameters (wavelength discretization~$N_\lambda$, Courant number~$S_\mathrm{c}$) and to the broader context of structure-preserving numerical methods, we aim to make this analysis both mathematically rigorous and immediately useful to practitioners.

\section{Geometry of the Problem and Electrodynamics Foundations}
\label{sect:Electrodynamics-Basics}

To facilitate comparison between the results derived in sections~\ref{sect:FDTD-Boundary-Conditions}--\ref{sect:Discussion} and their exact analytical counterparts, and to ensure the text is self-contained and readable, we begin with a concise review of key electrodynamic principles.
A more detailed treatment of these topics can be found in standard electrodynamics textbooks (e.g., \cite{LandauV8} or \cite{Jackson}).

\begin{condition}
\label{cond:Perfect-Media}
Throughout this work, a perfect medium refers to a lossless material where energy dissipation is entirely neglected (\mbox{$\varepsilon_\mathrm{r},\mu_\mathrm{r}\in\mathbb{R}$}).
Furthermore, we restrict our analysis to linear, homogeneous, isotropic, stationary media, treating~$\varepsilon_\mathrm{r}$ and~$\mu_\mathrm{r}$ as real constants~\cite{LandauV8,Jackson}.
Additionally, we focus exclusively on right-handed media~\cite{Veselago1967,Pendry2004}, where~$\varepsilon_\mathrm{r}$ and~$\mu_\mathrm{r}$ are simultaneously positive.
\end{condition}

\subsection{Problem Geometry}
\label{sect:ProblemGeometry}

Consider a plane-polarized electromagnetic wave incident normally upon a planar interface between two media.
We adopt a Cartesian coordinate system in which the interface lies at \mbox{$x=0$}.
The incident wave, with electric field~$\mathbf{E}^{(i)}$ and magnetic field~$\mathbf{H}^{(i)}$, propagates in the region of \mbox{$x<0$} along the $x$-axis (wavevector \mbox{$\mathbf{k}^{(i)}>0$}).
Because of the normal incidence of the plane-polarized wave, the problem reduces to one-dimensional dynamics in ($1+1$)-dimensional spacetime.
The $z$-axis is aligned with the direction of~$\mathbf{E}^{(i)}$, and the $y$-axis is fixed by the right-hand rule.

\begin{condition}
\label{cond:Currents-and-Charges-Absence}
External currents~$\mathbf{j}_\mathrm{ext}$ and free charges~$\rho_\mathrm{ext}$ are absent throughout the entire domain \mbox{$\mathcal{D}\subset\mathbb{R}^3$}, including the interface.
\end{condition}

\subsection{Governing Equations}
\label{sect:GoverningEquations}

Under Condition~\ref{cond:Perfect-Media} (lossless, linear, homogeneous, isotropic media) and Condition~\ref{cond:Currents-and-Charges-Absence}, Maxwell’s equations reduce to the following form:
\begin{equation}
\label{eq:Maxwells-Eqs}
  \varepsilon_0\varepsilon_\mathrm{r} \frac{\partial E_z}{\partial \mathrm{t}} = 
  \frac{\partial H_y}{\partial x}, \quad
  \mu_0\mu_\mathrm{r} \frac{\partial H_y}{\partial \mathrm{t}} = 
  \frac{\partial E_z}{\partial x}.
\end{equation}
Here, $\varepsilon_0$ and~$\mu_0$ are the vacuum permittivity and permeability, respectively, related to the speed of light~$c$ by:
\begin{equation}
\label{eq:Light-Speed}
  c = \frac{1}{\sqrt{\varepsilon_0\mu_0}}.
\end{equation}

The first equation in~\eqref{eq:Maxwells-Eqs} represents Amp\`{e}re’s law (with Maxwell’s displacement current), while the second encodes Faraday’s law.
A particular solution of~\eqref{eq:Maxwells-Eqs} is the harmonic plane wave:
\begin{equation}
\label{eq:Plane-Waves}
  E_z(x,\mathrm{t}) = E_+ \, e^{i(\omega \mathrm{t} - kx)} + E_- \, e^{i(\omega \mathrm{t} + kx)},\quad
  H_y(x,\mathrm{t}) = H_+ \, e^{i(\omega \mathrm{t} - kx)} + H_- \, e^{i(\omega \mathrm{t} + kx)},
\end{equation}
with the dispersion relation:
\begin{equation}
\label{eq:Dispersion-Law}
  k = \omega \frac{n_\mathrm{r}}{c},
\end{equation}
where $n_\mathrm{r}$ is the relative refractive index of the medium
\begin{equation}
\label{eq:Refraction-Index}
  n_\mathrm{r} \equiv \sqrt{\varepsilon_\mathrm{r}\mu_\mathrm{r}}.
\end{equation}

Here, $E_+$ and $E_-$ denote the amplitudes of waves propagating along and against the $x$-axis, respectively.
From~\eqref{eq:Maxwells-Eqs}\,--\,\eqref{eq:Dispersion-Law} , the magnetic field amplitudes~$H_\pm$ satisfy:
\begin{equation}
\label{eq:H-E-Relation}
  H_\pm = \mp \frac{E_\pm}{\eta},
\end{equation}
where~$\eta$ is the wave impedance of the medium (and~$\eta_0$---impedance of free space):
\begin{equation}
\label{eq:Impedance}
  \eta = \eta_0 \eta_\mathrm{r},\quad
  \eta_0 = \sqrt{\frac{\mu_0}{\varepsilon_0}},\quad
  \eta_\mathrm{r} = \sqrt{\frac{\mu_\mathrm{r}}{\varepsilon_\mathrm{r}}}
\end{equation}

The vectors $\mathbf{E}$, $\mathbf{H}$, and~$\mathbf{k}$ form a right-handed triad in right-handed media, which determines the relative signs in~\eqref{eq:H-E-Relation}.

\subsection{Reflection and Transmission at the Interface}
\label{sect:ReflectionAndTransmission}

\begin{figure}
\centering
\includegraphics[width=0.5\textwidth]{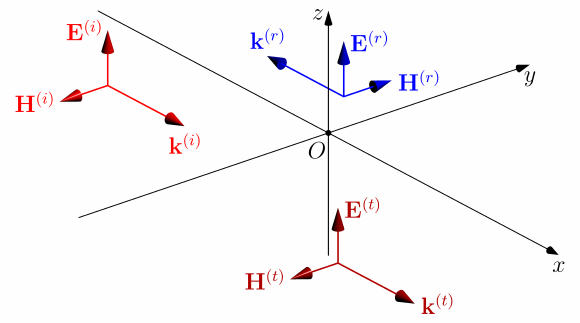}
\caption{Orientation of electromagnetic field components in incident, reflected and transmitted waves}
\label{fig:General-Geometry}
\end{figure}

At \mbox{$x=0$}, the interface separates Medium~1 (\mbox{$x<0$}: $\varepsilon_1$, $\mu_1$) from Medium~2 (\mbox{$x>0$}: $\varepsilon_2$, $\mu_2$).
The incident harmonic plane wave generates reflected ($\mathbf{E}^{(r)}$, $\mathbf{H}^{(r)}$) and transmitted ($\mathbf{E}^{(t)}$, $\mathbf{H}^{(t)}$) waves, as illustrated in Fig.~\ref{fig:General-Geometry}.
The fields are expressed as:
\begin{itemize}
\item Incident wave (\mbox{$x<0$}):
\begin{equation}
\label{eq:Incident-Fields}
  E_z^{(i)}(x,\mathrm{t}) = E_0 \, e^{i(\omega \mathrm{t} - k_1x)},\quad
  H_y^{(i)}(x,\mathrm{t}) = -\frac{E_0}{\eta_1} \, e^{i(\omega \mathrm{t} - k_1x)}.
\end{equation}
\item Reflected wave (\mbox{$x<0$}):
\begin{equation}
\label{eq:Reflected-Fields}
  E_z^{(r)}(x,\mathrm{t}) = r E_0 \, e^{i(\omega \mathrm{t} + k_1x)},\quad
  H_y^{(r)}(x,\mathrm{t}) = r \frac{E_0}{\eta_1} \, e^{i(\omega \mathrm{t} + k_1x)}.
\end{equation}
\item Transmitted wave (\mbox{$x>0$}):
\begin{equation}
\label{eq:Transmitted-Fields}
  E_z^{(t)}(x,\mathrm{t}) = t E_0 \, e^{i(\omega \mathrm{t} - k_2x)},\quad
  H_y^{(t)}(x,\mathrm{t}) = -t \frac{E_0}{\eta_2} \, e^{i(\omega \mathrm{t} - k_2x)}.
\end{equation}
\end{itemize}

Here, $r$ and~$t$ are Fresnel coefficients, defined by the relations
\begin{equation}
\label{eq:Fresnel-Coeffs-Def}
  r = \left.\frac{E^{(r)}}{E^{(i)}}\right|_{x=0}, \quad
  t = \left.\frac{E^{(t)}}{E^{(i)}}\right|_{x=0}.
\end{equation}

The coincidence of the angular frequencies~$\omega$ of the initial, reflected and transmitted waves is guaranteed by the linearity of the media considered within the Condition~\ref{cond:Perfect-Media}.

\subsection{Boundary Conditions and Exact Solution}
\label{sect:BoundaryConditions}

Continuity of the tangential field components at \mbox{$x=0$} (Condition~\ref{cond:Currents-and-Charges-Absence}) yields:
\begin{equation}
\label{eq:Boundary-Conditions}
  E_z^{(i)}(0,\mathrm{t}) + E_z^{(r)}(0,\mathrm{t}) = E_z^{(t)}(0,\mathrm{t}),\quad
  H_y^{(i)}(0,\mathrm{t}) + H_y^{(r)}(0,\mathrm{t}) = H_y^{(t)}(0,\mathrm{t}).
\end{equation}

Substituting~\eqref{eq:Incident-Fields}\,--\,\eqref{eq:Transmitted-Fields} into~\eqref{eq:Boundary-Conditions} and solving for~$t$ and~$r$, we obtain:
\begin{equation}
\label{eq:t-r-Exact}
  t = \frac{2\eta_2}{\eta_2 + \eta_1},\quad
  r = \frac{\eta_2 - \eta_1}{\eta_2 + \eta_1}.
\end{equation}

The Fresnel coefficients~\eqref{eq:t-r-Exact} relate the amplitudes of the incident, reflected and refracted waves in the case of normal incidence of a harmonic plane wave at the interface of two media with impedances~$\eta_1$ and~$\eta_2$.
In practice it is necessary to estimate not only the amplitudes of the corresponding waves, but the average values of their energy fluxes.
This problem is solved with the reflection and transmission coefficients, related to the Fresnel coefficients in the considered case of perfect media (Condition~\ref{cond:Perfect-Media}) as follows:
\begin{equation}
\label{eq:R-T-Exact}
  R = r^2,\quad
  T = \frac{\eta_1}{\eta_2} t^2.
\end{equation}

From~\eqref{eq:t-r-Exact} and~\eqref{eq:R-T-Exact}, it is easy to see that the equality \mbox{$R+T=1$} is satisfied, which is a one form of the law of conservation of electromagnetic energy.
The exact solutions~\eqref{eq:t-r-Exact} and~\eqref{eq:R-T-Exact} will serve as a benchmark in \S\ref{sect:Discussion} for evaluating the accuracy of FDTD simulations.

\section{Space-time discretization and finite-differences approach}
\label{sect:FDTD-Basics}

Our next objective is to briefly outline the FDTD method, on which all subsequent results are based.
The FDTD method, first proposed by Kane Yee in his pioneering work~\cite{Yee1966}, is essentially a technique for approximating partial derivatives of continuous functions by finite differences of corresponding quantities
\begin{gather}
  \label{eq:Dpsi-Dx}
  \frac{\partial \psi(x,\mathrm{t})}{\partial x}\Big|_{(x_0,\mathrm{t}_0)} =
  \frac{\psi(x_0 + \Delta_x,\mathrm{t}_0) - \psi(x_0,\mathrm{t}_0)}{\Delta_x} + o(1),\\
  \label{eq:Dpsi-Dt}
  \frac{\partial \psi(x,\mathrm{t})}{\partial \mathrm{t}}\Big|_{(x_0,\mathrm{t}_0)} =
  \frac{\psi(x_0,\mathrm{t}_0 + \Delta_t) - \psi(x_0,\mathrm{t}_0)}{\Delta_t} + o(1).
\end{gather}
Here and below, $\psi$ denotes an arbitrary component of the electromagnetic field: $E_z$ or~$H_y$.
The above relations approximate the first-order derivatives of~$\psi$ using the forward difference to first order in the step size between adjacent values of~$\psi$.

\subsection{Yee Algorithm}
\label{sect:YeeAlgorithm}

The adoption of approximation~\eqref{eq:Dpsi-Dx}, \eqref{eq:Dpsi-Dt} naturally leads to the necessity of transition from continuous functions of space-time coordinates~$\psi(x,t)$ to their discrete analogues~$\psi^q[m]$.
For this purpose, a time-space grid with constant time~$\Delta_t$ and spatial~$\Delta_x$ steps is introduced, on the nodes of which, using a formal transition
\begin{equation}
  \psi(x,\mathrm{t}) = \psi(m\Delta_x, q\Delta_t) = \psi^q[m],
\end{equation}
the values of the nodal function~$\psi^q[m]$ are determined.
As can be seen from the above equality, the index $m=0,1,\dots,M-1$ denotes the spatial position of the grid node ($M$ is number of spatial nodes of the Yee grid in the $x$-axis direction of the domain~$\mathcal{D}$), while $q=0,1,\dots,T-1$ indexes moment of time ($T$ is the number of time steps that regulate the total duration of the simulation).

\begin{figure}
\centering
\includegraphics[width=0.5\textwidth]{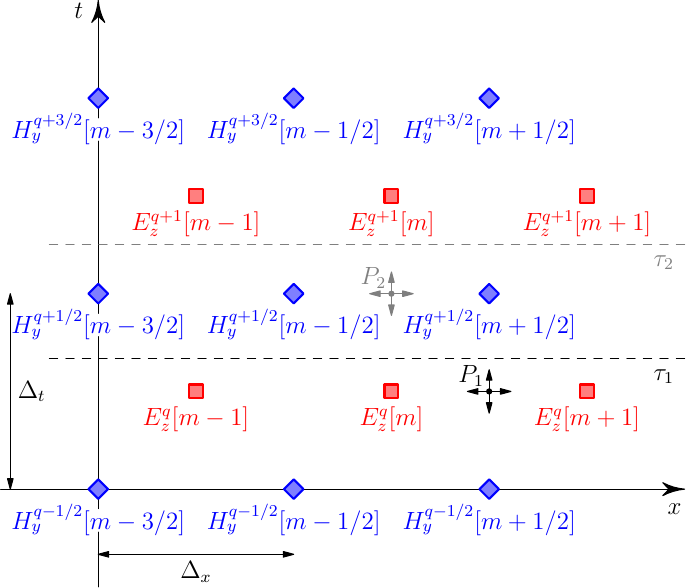}
\caption{Yee grid}
\label{fig:Yee-Grid}
\end{figure}

The main feature of the FDTD method proposed by Yee is a staggered time-space discretization in which the grids for the electric~$E_z$ and magnetic~$H_y$ components of electromagnetic field are shifted relative to each other by half a step in each dimension.
This picture is shown in Fig.~\ref{fig:Yee-Grid}.
From this it is evident that the total number of nodes in the Yee grid is \mbox{$2M\times2T$}, which are not completely independent of each other, since half of them are associated with the electric component of the field, and half---with the magnetic part.
Thus, the total spatial extent of the considered region of space~$\mathcal{D}$ is determined precisely by the number~$M$, and the duration of the simulation is fixed by the number~$T$.

Let us now refine the approximation~\eqref{eq:Dpsi-Dx}, \eqref{eq:Dpsi-Dt} using central differences at the spacetime point with coordinates~$(x_0,\mathrm{t}_0)$ shown in Fig.~\ref{fig:Yee-Grid} at position~$P_1$ with discrete indices~$m+1/2$ and~$q$.
The Taylor series expansion of the functions~$\psi$ at the points $(x_0\pm\Delta_x/2,\mathrm{t}_0)$ and~$(x_0,\mathrm{t}_0\pm\Delta_t/2)$ with the accuracy up to the second order of smallness allows to write down
\begin{gather}
  \label{eq:Dpsi-Dx-2}
  \frac{\partial \psi(x,\mathrm{t})}{\partial x}\Big|_{(x_0,\mathrm{t}_0)} =
  \frac{\psi(x_0 + \Delta_x/2,\mathrm{t}_0) - \psi(x_0 - \Delta_x/2,\mathrm{t}_0)}{\Delta_x} + o(\Delta_x),\\
  \label{eq:Dpsi-Dt-2}
  \frac{\partial \psi(x,\mathrm{t})}{\partial \mathrm{t}}\Big|_{(x_0,\mathrm{t}_0)} =
  \frac{\psi(x_0,\mathrm{t}_0 + \Delta_t/2) - \psi(x_0,\mathrm{t}_0 - \Delta_t/2)}{\Delta_t} + o(\Delta_t).
\end{gather}

Using approximations~\eqref{eq:Dpsi-Dx-2}, \eqref{eq:Dpsi-Dt-2} for point~$P_1$ yields the FDTD-analogue of the second equation in~\eqref{eq:Maxwells-Eqs} in the following form
\begin{equation}\label{eq:FDTD-Faraday-Law}
  H_y^{q+\frac{1}{2}}\left[m+\frac{1}{2}\right] =
  H_y^{q-\frac{1}{2}}\left[m+\frac{1}{2}\right] +
  \frac{S_\mathrm{c}}{\eta_0\mu_\mathrm{r}}\left(
  E_z^q[m+1] - E_z^q[m]\right).
\end{equation}

Here, when writing the finite-difference analogue of Faraday's law, the notation for the vacuum impedance~\eqref{eq:Impedance}, the definition of the Courant number
\begin{equation}
\label{eq:Courant-Number}
  S_\mathrm{c} \equiv \frac{c\Delta_t}{\Delta_x},
\end{equation}
and the relationship between the material parameters of the vacuum with the speed of light~\eqref{eq:Light-Speed} were used.
In addition, Equation~\eqref{eq:FDTD-Faraday-Law} is presented in such a form that it allows one to express the ``future'' values of~$H_y^{q+1/2}$, $\forall m$, if ``past'' values of~$E_z^q$ and~$H_y^{q-1/2}$ are known.
The dashed line~$\tau_1$ in Fig.~\ref{fig:Yee-Grid} marking the boundary dividing the ``past'' and ``future'' regions for node~$P_1$.

By carrying out a similar approximation procedure~\eqref{eq:Dpsi-Dx-2}, \eqref{eq:Dpsi-Dt-2} for the time-space grid node~$P_2$, from the first equation in~\eqref{eq:Maxwells-Eqs} one can obtain the FDTD-analogue of the Amp\`{e}re-Maxwell law equation
\begin{equation}\label{eq:FDTD-Ampere-Law}
  E_z^{q+1}[m] = E_z^{q}[m] +
  \frac{S_\mathrm{c}\eta_0}{\varepsilon_\mathrm{r}}\left(
  H_y^{q+\frac{1}{2}}\left[m+\frac{1}{2}\right] -
  H_y^{q+\frac{1}{2}}\left[m-\frac{1}{2}\right]\right).
\end{equation}

The~\eqref{eq:FDTD-Faraday-Law} and~\eqref{eq:FDTD-Ampere-Law} obtained above represent a coupled system of equations that allows the implementation of a recurrent scheme for determining all components of the electromagnetic field using known initial conditions.
For this reason, that in the literature on the FDTD method these equations are known as the Update Equations~\cite{Schneider,Inan}.

The Yee scheme’s staggered grid (Fig.~\ref{fig:Yee-Grid}) ensures second-order accuracy in space and time.
Update equations~\eqref{eq:FDTD-Faraday-Law}, \eqref{eq:FDTD-Ampere-Law} decouple the electric and magnetic fields, preserving the symplectic structure of Maxwell’s equations~\eqref{eq:Maxwells-Eqs}.

Our recent work~\cite{Makarov2024a} has demonstrated that setting the Yee grid parameters as
\begin{equation}
\label{eq:Magic-Courant-Number}
  S_\mathrm{c} = n_\mathrm{r},
\end{equation}
is optimal for modeling the electromagnetic properties of a homogeneous, non-dispersive medium.

This choice significantly reduces computational errors arising from numerical dispersion when simulating wave propagation in a medium with relative permittivity~$\varepsilon_\mathrm{r}$ and permeability~$\mu_\mathrm{r}$, modeled as a discrete nodal chain.

\subsection{Generation of the Incident Wave Field}
\label{sect:IncidentWaveGeneration}

From Eqs.~\eqref{eq:FDTD-Faraday-Law} and~\eqref{eq:FDTD-Ampere-Law}, it follows that in the framework of the Condition~\ref{cond:Currents-and-Charges-Absence} accepted by us, nontrivial dynamics of the electromagnetic field is possible only in two cases~\cite{Makarov2024b}:
\begin{enumerate}
  \item Static mode. It is realized when setting a non-zero initial configuration of the electromagnetic field \mbox{$(\psi^0\cup\psi^{1/2})\not\equiv0$}, \mbox{$\forall m\in[0,M)$}.
  \item Dynamic mode. The electromagnetic field of the incident wave must be correctly described at each time step of the simulation~$q$.
\end{enumerate}

\begin{remark}
The first of the above variants is equivalent to the introduction of such external electric charges \mbox{$\rho_\mathrm{s}\neq0$}, that act as sources setting its initial state.
Therefore, the external free charges~$\rho_\mathrm{ext}$ in the Condition~\ref{cond:Currents-and-Charges-Absence} are understood as all charges except for the given charges~$\rho_\mathrm{s}$, which we have specifically mentioned here.
The second variant is equivalent to the introduction of external current sources~$\mathbf{j}_\mathrm{s}$ controlling the behavior of the electromagnetic field in the dynamics.
Here it should also be clarified that the external currents~$\mathbf{j}_\mathrm{ext}$ in the Condition~\ref{cond:Currents-and-Charges-Absence} are any currents additional to~\mbox{$\mathbf{j}_s$}.
\end{remark}

In this paper, we will define the field of the incident wave in the dynamic mode, assuming that its source is the current~$j_z^q[s]$ generated by the transmitter~$J$---a point directivity antenna.
It is assumed that the transmitter~$J$ is located on a node of the Yee grid with spatial index \mbox{$m=s$}.
By the directionality of the source~$J$ we mean that its action forms an initial wave~\eqref{eq:Incident-Fields} propagating in the direction of \mbox{$m\geqslant s$}.

In the presence of inhomogeneities (in particular, the interfaces between different media, which is the subject of the present study), when implementing the iterative algorithm~\eqref{eq:FDTD-Faraday-Law}, \eqref{eq:FDTD-Ampere-Law} in the area \mbox{$m\geqslant s$} of the Yee grid with the flow of time \mbox{$q\geqslant0$} the total field, which is the sum of the incident and reflected waves, is formed stepwise.
In the spatial region of the Yee grid \mbox{$m<s$}, only the so-called scattered field (created in the case under consideration by the reflected wave~\eqref{eq:Reflected-Fields}) can penetrate.
For this reason, the dynamical mode of setting the nontrivial dynamics of the electromagnetic field that we describe here is called Total Field/Scattered Field (TFSF) in the FDTD literature~\cite{Schneider,Inan,Langtangen}.

Without going further into more detailed explanations (which can be found, for example, in the excellent book~\cite{Schneider}), we will only point out that the action of current~$j_z^q[s]$, consistent with the choice of the direction of the initial wave~\eqref{eq:Incident-Fields}, in the TFSF formalism is reduced to taking into account summands additional to the right-hand sides~\eqref{eq:FDTD-Faraday-Law} and~\eqref{eq:FDTD-Ampere-Law}, for two nodes of the grid:
\begin{equation}\label{eq:Source-Hy}
  H_y^{q+\frac12}\!\left[s-\frac12\right] =
  -\frac{S_\mathrm{c}}{\eta_0\mu_\mathrm{r}} f^{q}\left[0\right],
\end{equation}
\begin{equation}\label{eq:Source-Ez}
  E_z^{q+1}[s] =
  \frac{S_\mathrm{c}}{\sqrt{\varepsilon_\mathrm{r}\mu_\mathrm{r}}}
  f^{q+\frac12}\!\left[-\frac12\right].
\end{equation}
Here~$f^q[m]$ is the finite-difference analog of the continuous electric component~$E_z^{(i)}(x,\mathrm{t})$ of the incident electromagnetic field~\eqref{eq:Incident-Fields}.

\subsection{FDTD Analogs of Incident, Reflected and Transmitted Wave Fields}
\label{sect:FDTD-Fields}

For future reference, we also list here the FDTD analogs of the fields~\eqref{eq:Incident-Fields}\,--\,\eqref{eq:Transmitted-Fields}.
Here and everywhere else in the following, we assume that the interface between Medium~1 and Medium~2 (plane \mbox{$x=0$}) is the Yee grid node with the number \mbox{$m=b$} (all necessary details will be set out below in section~\ref{sect:FDTD-Boundary-Conditions}).
\begin{itemize}
\item Incident wave (\mbox{$s\leqslant m\leqslant b$}):
\begin{equation}\label{eq:FDTD-Incident-Fields}
  f^q[m] \equiv E_z^{(i)\,q}[m] =
  E_0 \, e^{i(\omega q\Delta_t - \widetilde{k}_1 (m-b)\Delta_x)},\quad
  H_y^{(i)\,q}[m] =
  -\frac{E_0}{\eta_1} \, e^{i(\omega q\Delta_t - \widetilde{k}_1 (m-b)\Delta_x)}.
\end{equation}
\item Reflected wave (\mbox{$m\leqslant b$}):
\begin{equation}\label{eq:FDTD-Reflected-Fields}
  E_z^{(r)\,q}[m] =
  \widetilde{r} E_0 \, e^{i(\omega q\Delta_t + \widetilde{k}_1 (m-b)\Delta_x)},\quad
  H_y^{(r)\,q}[m] =
  \widetilde{r} \, \frac{E_0}{\eta_1} \, e^{i(\omega q\Delta_t + \widetilde{k}_1 (m-b)\Delta_x)}.
\end{equation}
\item Transmitted wave (\mbox{$m\geqslant b$}):
\begin{equation}\label{eq:FDTD-Transmitted-Fields}
  E_z^{(t)\,q}[m] =
  \widetilde{t} E_0 \, e^{i(\omega q\Delta_t - \widetilde{k}_2 (m-b)\Delta_x)},\quad
  H_y^{(t)\,q}[m] =
  -\widetilde{t} \, \frac{E_0}{\eta_2} \, e^{i(\omega q\Delta_t - \widetilde{k}_2 (m-b)\Delta_x)}.
\end{equation}
\end{itemize}

\begin{remark}
Here and below, the tilde sign written above the symbol indicates that this quantity is a FDTD analogue of the corresponding continuous quantity, which in the general case does not coincide with it (for example, for FDTD-wavenumbers~$\widetilde{k}$ see the section~\ref{sect:FDTD-Dispersion}, and the obtaining of the quantities~$\widetilde{r}$ and~$\widetilde{t}$ is essentially one of the main goals of the work and it is described in detail below in sections~\ref{sect:FDTD-Boundary-Conditions}--\ref{sect:Discussion}).
\end{remark}

\begin{remark}
We note separately that the magnitude~$E_0$ and angular frequency~$\omega$ of a harmonic incident wave is assumed to be a given quantity for which there is no difference between the continuous and discretized versions.
At the same time, the FDTD-impedance matches with its continuous counterpart, the proof of which can be found in~\cite[7.5]{Schneider}.
Because of that, these quantities are written in~\eqref{eq:FDTD-Incident-Fields}\,--\,\eqref{eq:FDTD-Transmitted-Fields} and below without tilde sign.
\end{remark}

\subsection{FDTD Dispersion}
\label{sect:FDTD-Dispersion}

The last thing we need to discuss before moving directly to achieving the goal of this study is the consideration of numerical dispersion, which inevitably takes place in FDTD simulations (even in the case of non-dispersive and homogeneous materials).
However, since this issue is studied in detail in the article~\cite{Makarov2024a} (also, see \cite[7.4]{Schneider}), here we will reproduce only the main result, which we will need later.

\begin{lemma}[see Statement\,1 in~\cite{Makarov2024a}]\label{lemma:FDTD-Dispersion}
The dispersion relation for the Yee grid can be represented in the form
\begin{equation}\label{eq:Yee-Dispersion-Eq}
  \sin \left( \frac{\widetilde{k}\Delta_x}{2} \right) =
  \frac{n_\mathrm{r}}{S_\mathrm{c}}
  \sin \left( \frac{\omega\Delta_t}{2} \right)
\end{equation}
\end{lemma}

\begin{proof}
Is given in the paper~\cite{Makarov2024a} and is not reproduced here. Expression~\eqref{eq:Yee-Dispersion-Eq} is written using the notation~\eqref{eq:Courant-Number}, \eqref{eq:Refraction-Index} and~\eqref{eq:Light-Speed} and rearranged as adopted in this paper.
\end{proof}

\begin{remark}
\label{note:Dispersion}
The dispersion relation~\eqref{eq:Yee-Dispersion-Eq} of the FDTD method in the general case (when the equality~\eqref{eq:Magic-Courant-Number} is not satisfied) differs significantly from its continuous counterpart, which has the form~\eqref{eq:Dispersion-Law}.
However, we note that this difference becomes vanishingly small at sufficiently small discretization of spacetime.
Indeed, by keeping only the first order of magnitude summands in the Taylor series of the sine expansion at~$\Delta_t$ and \mbox{$\Delta_x\rightarrow0$}, it is easy to obtain from~\eqref{eq:Yee-Dispersion-Eq}
\begin{equation}\label{eq:Yee-Dispersion-DtDx0}
  \widetilde{k} = \omega \frac{n_\mathrm{r}}{c}.
\end{equation}

We emphasize, however, that~\eqref{eq:Yee-Dispersion-DtDx0} is valid only in the limit of infinitesimal discretization \mbox{$\Delta_t,\Delta_x\rightarrow0$}.
In the general case of finite discretization of spacetime, the expression~\eqref{eq:Yee-Dispersion-Eq} must be used to determine~$\widetilde{k}$.
\end{remark}

The interplay between numerical dispersion and stability constraints remains an active topic in the analysis of time-domain methods for hyperbolic systems; see, e.g., recent studies on energy-stable discretizations \cite{Zhang2024}, high-order compact schemes \cite{Kong2025} and quantitative error estimates for an ADI scheme \cite{Liu2025} in \emph{Applied Numerical Mathematics}.

\section{FDTD boundary conditions}
\label{sect:FDTD-Boundary-Conditions}

\subsection{Necessary definitions and lemmas}
\label{sect:FDTD-DefinitionsAndLemmas}

To obtain the expressions describing the boundary conditions for the Yee grid, we return to the original discrete forms of the Amp\`{e}re~\eqref{eq:FDTD-Ampere-Law} and Faraday~\eqref{eq:FDTD-Faraday-Law} laws, which in a more compact and convenient form can be written using the shift operators in the space-time grid~$\widehat{S}_x^\chi$ and~$\widehat{S}_t^\tau$, whose action we define as follows.

\begin{definition}
\label{def:S}
Let~$\chi$ and~$\tau$---such numbers that can be represented in the form \mbox{$\chi,\tau=p/2,\,\forall p\in\mathbb{Z}$}. Then
\begin{gather}
  \label{eq:Sx-Def}
  \widehat{S}_x^\chi \psi^q[m] = \psi^q\left[m+\chi\right],\\
  \label{eq:St-Def}
  \widehat{S}_t^\tau \psi^q[m] = \psi^{q+\tau}[m].
\end{gather}
\end{definition}

\begin{example}
Let us write here explicitly the action of the shift operators~\eqref{eq:Sx-Def},\eqref{eq:St-Def} on an arbitrary component~$\psi$ of the electromagnetic fields~\eqref{eq:FDTD-Incident-Fields}\,--\,\eqref{eq:FDTD-Transmitted-Fields} only with parameters \mbox{$\chi,\tau=\pm1/2$}
\begin{equation}\label{eq:Sx-St-Plane-Wave}
  \widehat{S}_t^{\pm\frac12} \psi =
  e^{\pm i\omega\Delta_t/2} \psi,\quad
  \widehat{S}_x^{\pm\frac12} \psi^{(i, t)} =
  e^{\mp i\widetilde{k}_{(1,2)}\Delta_x/2} \psi^{(i, t)},\quad
  \widehat{S}_x^{\pm\frac12} \psi^{(r)} =
  e^{\pm i\widetilde{k}_1\Delta_x/2} \psi^{(r)},
\end{equation}
\end{example}

It can be easily established that for the operators~\eqref{eq:Sx-Def}, \eqref{eq:St-Def} the following property is true
\begin{equation}
\label{eq:S-Reciprocity}
  \widehat{S}_x^\chi \widehat{S}_x^{-\chi} =
  \widehat{S}_x^{-\chi} \widehat{S}_x^\chi =
  \widehat{S}_x^0 \equiv \widehat{I},\quad
  \widehat{S}_t^\tau \widehat{S}_t^{-\tau} =
  \widehat{S}_t^{-\tau} \widehat{S}_t^\tau =
  \widehat{S}_t^0 \equiv \widehat{I},
\end{equation}
where~$\widehat{I}$ is a unit operator: \mbox{$\widehat{I}\psi^q[m]=\psi^q[m]$}.

\begin{definition}
\label{def:D}
Let us also define for convenience the finite-difference operators~$\widetilde{\partial}_i$ according to
\begin{equation}
\label{eq:D-Operator-Def}
  \widetilde{\partial}_i =
  \frac{\widehat{S}_i^{\frac12} - \widehat{S}_i^{-\frac12}}{\Delta_i},
\end{equation}
where~$i$ is~$x$ or~$t$.
\end{definition}

\begin{lemma}[{see~\cite[7.3, 7.4]{Schneider}}]\label{lemma:Yee-Form-of-Maxwells-Eqs}
With~\eqref{eq:D-Operator-Def}, the FDTD analog of Amp\`{e}re's law under Conditions~\ref{cond:Perfect-Media} and~\ref{cond:Currents-and-Charges-Absence} can be written in the Yee form
\begin{equation}\label{eq:Ampere-Law-Yee}
  \varepsilon_0\varepsilon_\mathrm{r}
  \widetilde{\partial}_t E_z^q[m] =
  \widetilde{\partial}_x H_y^q[m].
\end{equation}
Similarly, the discrete analog of Faraday's law in the form of Yee takes the following form
\begin{equation}\label{eq:Faraday-Law-Yee}
  \mu_0\mu_\mathrm{r}
  \widetilde{\partial}_t H_y^q[m] =
  \widetilde{\partial}_x E_z^q[m].
\end{equation}
\end{lemma}

\begin{proof}
The equality~\eqref{eq:Ampere-Law-Yee} follows directly from the first equation of the system~\eqref{eq:Maxwells-Eqs} after applying the approximation~\eqref{eq:Dpsi-Dx-2}, \eqref{eq:Dpsi-Dt-2} for the node of the Yee grid with coordinates~$m$, $q+1/2$, labeled in Fig.~\ref{fig:Yee-Grid} as point~$P_2$, with the use of Definitions~\ref{def:S} and~\ref{def:D}.
\begin{equation*}
  \varepsilon_0\varepsilon_\mathrm{r} \widehat{S}_t^{\frac12}
  \widetilde{\partial}_t E_z^q[m] =
  \widehat{S}_t^{\frac12} \widetilde{\partial}_x H_y^q[m].
\end{equation*}
On the next step, let us act on the left and right parts of above equality with the operator~$\widehat{S}_t^{-\frac12}$ and use its reciprocity property~\eqref{eq:S-Reciprocity}.
These steps lead us to~\eqref{eq:Ampere-Law-Yee}.

The same actions performed for the second equation in~\eqref{eq:Maxwells-Eqs} at the point~$P_1$ of Fig.~\ref{fig:Yee-Grid} lead to the~\eqref{eq:Faraday-Law-Yee}.
\end{proof}

\begin{lemma}
\label{lemma:D-Actions}
The action of the operators~\eqref{eq:D-Operator-Def} on an arbitrary component~$\psi$ of harmonic plane waves~\eqref{eq:FDTD-Incident-Fields}\,--\,\eqref{eq:FDTD-Transmitted-Fields} is
\begin{equation}\label{eq:Fin-Diff-Operator-Dt}
  \widetilde{\partial}_t \psi = i \Omega \psi,\quad
  \Omega \equiv \frac{2}{\Delta_t} \sin \left( \frac{\omega\Delta_t}{2} \right),
\end{equation}
\begin{equation}\label{eq:Fin-Diff-Operator-Dx}
  \widetilde{\partial}_x \psi^{(i,t)} = -i K_{(1,2)} \psi^{(i,t)},\quad
  \widetilde{\partial}_x \psi^{(r)} = i K_1 \psi^{(r)},\quad
  K_{(1,2)} \equiv \frac{2}{\Delta_x} \sin \left( \frac{\widetilde{k}_{(1,2)}\Delta_x}{2} \right).
\end{equation}
\end{lemma}

\begin{proof}
Verified by substituting~\eqref{eq:FDTD-Incident-Fields}\,--\,\eqref{eq:FDTD-Transmitted-Fields} into~\eqref{eq:D-Operator-Def} taking into account~\eqref{eq:Sx-St-Plane-Wave}.
\end{proof}

Now, after discussing all necessary auxiliary details, we proceed to the direct derivation of the boundary conditions for the interface of two media.
Under the assumption of Condition~\ref{cond:Perfect-Media}, it is clear that the relative permittivity~$\varepsilon_\mathrm{r}$ and permeability~$\mu_\mathrm{r}$ are piecewise constant functions of the Yee grid spatial node: \mbox{$\varepsilon_\mathrm{r}=\varepsilon_\mathrm{r}[m]$} and \mbox{$\mu_\mathrm{r}=\mu_\mathrm{r}[m]$}.
Moreover, it is evident from~\eqref{eq:Ampere-Law-Yee} and~\eqref{eq:Faraday-Law-Yee} that the function~$\varepsilon_\mathrm{r}[m]$ must be linked to the electric field~$E^q_z[m]$, and the nodes of the function~$\mu_\mathrm{r}[m]$ need to be connected to the nodes of the~$H^q_y[m]$ function.

Let us now examine in detail the cases of dielectrics and magnetics separately.

\subsection{Dielectrics}
\label{sect:DielectricsBoundary}

\begin{figure}
\centering
\includegraphics[width=0.75\textwidth]{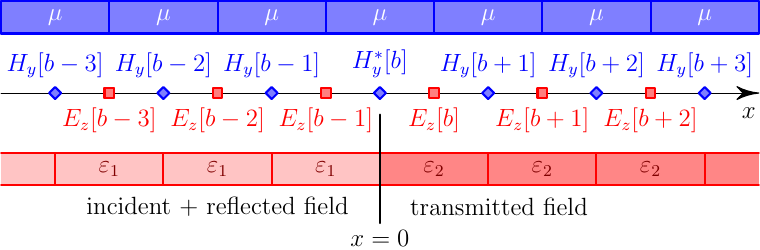}
\caption{FDTD-model of a planar interface between two dielectrics}
\label{fig:Dielectrics-Interface}
\end{figure}

Consider the interface of two dielectrics with relative permittivities~$\varepsilon_1$ and~$\varepsilon_2$.
Let the relative permeability of both media be the same and equal to~$\mu$ (for generality of the problem formulation $\mu$ is not necessarily equal to~$1$).
Geometry of this model is shown in Fig.~\ref{fig:Dielectrics-Interface}.
On that scheme we assume that Medium~1 extends over~$b$ nodes, such that
\begin{equation*}
  \varepsilon_\mathrm{r}[m] = \begin{cases}
  \varepsilon_1, m\in[0,b)\\
  \varepsilon_2, m\in[b,M)
  \end{cases},\quad
  \mu_\mathrm{r}[m] \equiv \mu = \const.
\end{equation*}

\begin{remark}
Additionally, we emphasize that under our assumptions the source~$J$ of the incident wave is located in the first medium, so that \mbox{$s<b$}.
\end{remark}

\begin{remark}
It should also be noted that the designations of the magnetic nodes in Fig.~\ref{fig:Dielectrics-Interface} and after this are shifted by a factor of~$1/2$.
This is done only for convenience, as the meaning remains the same as in Fig.~\ref{fig:Yee-Grid}.
Thus, here and further, the nodes~$E^q_z[p]$ and~$H^q_y[p]$ with the same spatial indexes~$p$ are actually shifted in space relative to each other by~$\Delta_x/2$.
\end{remark}

Under the conditions we have adopted here, we can see from Fig.~\ref{fig:Dielectrics-Interface} that the interface between Medium~1 and Medium~2 (plane \mbox{$x=0$}) in the discrete space model of the FDTD method should be chosen as the magnetic field node $H_y[b]$, marked with an asterisk in this figure.
For this node we can easily obtain the first FDTD boundary condition for dielectrics---the continuity of the tangential component of the magnetic field (second equation in~\eqref{eq:Boundary-Conditions}).
After substitution~\eqref{eq:FDTD-Incident-Fields}\,--\,\eqref{eq:FDTD-Transmitted-Fields} into~\eqref{eq:Boundary-Conditions} and canceling the common multiplier~$E_0 e^{i\omega q\Delta_t}$, this leads to the equality
\begin{equation}
\label{eq:Continity-Condition-1}
  -\frac1{\eta_1} + \frac{\widetilde{r}}{\eta_1} =
  -\frac{\widetilde{t}}{\eta_2}.
\end{equation}

The second FDTD boundary condition is more difficult to obtain.
This follows from the principal absence of exact discrete analog of the first equation of the system~\eqref{eq:Boundary-Conditions} due to staggered placement of Yee's grid nodes.
This statement is clearly demonstrated by Fig.~\ref{fig:Dielectrics-Interface}, from which it is obvious that within the framework of the model we have adopted, there is no electric node of the Yee grid lying exactly at the interface between the two media.
Thus, the nodes of the Yee grid~$E_z[b-1]$ and~$E_z[b]$, nearest to the interface, lie strictly in the Medium~1 and Medium~2, respectively, and therefore cannot accurately link together the incident, reflected, and transmitted waves.
The formulated problem is solved by the following theorem.

\begin{theorem}
\label{thrm:1}
The FDTD-analog of the boundary condition relating the tangential components of the electric field of a plane monochromatic wave at the interface of two dielectrics within the model presented in Fig.~\ref{fig:Dielectrics-Interface}, is the following equation
\begin{equation}\label{eq:Continity-Condition-2}
  e^{i\widetilde{k}_1\Delta_x/2} +
  \widetilde{r}e^{-i\widetilde{k}_1\Delta_x/2} =
  \left(e^{-i\widetilde{k}_2\Delta_x/2} +
  i\frac{\Omega\mu_0\mu\Delta_x}{\eta_2}\right) \! \widetilde{t}.
\end{equation}
\end{theorem}

\begin{proof}
First of all, in this proof we rely on the result of Lemma~\ref{lemma:Yee-Form-of-Maxwells-Eqs}.
Namely, we write the FDTD-analog of Faraday's law in the form~\eqref{eq:Faraday-Law-Yee} for the boundary node~$H^*_y[b]$, since by using this equality we can approximately relate together components of electric field in Medium~1 and Medium~2.
\begin{equation}
\label{eq:Theorem1-Proof-Step-1}
  \mu_0 \mu \widetilde{\partial}_t H_y^q[b] =
  \widetilde{\partial}_x E_z^q[b].
\end{equation}

Then, we use the definition~\eqref{eq:D-Operator-Def} of the operator~$\widetilde{\partial}_x$ and the result~\eqref{eq:Fin-Diff-Operator-Dt} of Lemma~\ref{lemma:D-Actions}.
Thus, we obtain the following relations
\begin{equation}
\label{eq:Theorem1-Proof-Step-2}
  i \Omega \mu_0 \mu H_y^q[b] =
  \frac1{\Delta_x}
  \left(\widehat{S}_x^{\frac12} - \widehat{S}_x^{-\frac12} \right) E_z^q[b].
\end{equation}

Let us now consider separately the action of the shift operators~$\widehat{S}^{\pm\frac12}_x$ in the right-hand side of the equality~\eqref{eq:Theorem1-Proof-Step-2} on the electric field~$E_z^q[b]$.
Thus,
\begin{equation}
\label{eq:Theorem1-Proof-Step-3}
  \widehat{S}_x^{\frac12} E_z^q[b] = e^{-\widetilde{k}_2\Delta_x/2} E_z^{(t)\,q}[b] =
  e^{-\widetilde{k}_2\Delta_x/2} \, \widetilde{t} \, E_0 \, e^{i\omega q\Delta_t},
\end{equation}
\begin{equation}
\label{eq:Theorem1-Proof-Step-4}
  \widehat{S}_x^{-\frac12} E_z^q[b] =
  e^{\widetilde{k}_1\Delta_x/2} E_z^{(i)\,q}[b] +
  e^{-\widetilde{k}_1\Delta_x/2} E_z^{(r)\,q}[b] =
  \left(e^{\widetilde{k}_1\Delta_x/2} + \widetilde{r} \, e^{-\widetilde{k}_1\Delta_x/2}\right)
  E_0 \, e^{i\omega q\Delta_t}.
\end{equation}
In writing these equations, it is taken into account that~$\widehat{S}_x^{\frac12}E_z^q[b]$ is the transmitted wave~\eqref{eq:FDTD-Transmitted-Fields}, and~$\widehat{S}_x^{-\frac12} E_z^q[b]$ is the sum of the incident~\eqref{eq:FDTD-Incident-Fields} and reflected~\eqref{eq:FDTD-Reflected-Fields} waves, and the facts~\eqref{eq:Sx-St-Plane-Wave} about the effect of the~$\widehat{S}^{\pm\frac12}_x$ operators, acting on the components of these waves are used.

Note now that the left part of the equality~\eqref{eq:Theorem1-Proof-Step-2} by virtue of the boundary condition~\eqref{eq:Continity-Condition-1} can be interpreted either as the sum of the incident and reflected waves or as a transmitted wave.
Taking into account the second of these options and using the expressions for the magnetic component of the electromagnetic field of the transmitted wave~\eqref{eq:FDTD-Transmitted-Fields}, we write the left part of~\eqref{eq:Theorem1-Proof-Step-2} in the form of
\begin{equation}
\label{eq:Theorem1-Proof-Step-5}
  i \Omega \mu_0 \mu H_y^q[b] =
  -i \Omega \mu_0 \mu \, \widetilde{t} \, \frac{E_0}{\eta_2} \, e^{i\omega q\Delta_t}.
\end{equation}

Substituting now~\eqref{eq:Theorem1-Proof-Step-3}\,--\,\eqref{eq:Theorem1-Proof-Step-5} into~\eqref{eq:Theorem1-Proof-Step-2}, reducing the common multiplier~$E_0e^{i\omega q\Delta_t}$, multiplying by~$\Delta_x$ and grouping the terms with coefficient~$\widetilde{t}$ in the right-hand side of the equality, we finally arrive at the expression~\eqref{eq:Continity-Condition-2}, which completes the proof.
\end{proof}

\subsection{Magnetics}
\label{sect:MagneticsBoundary}

\begin{figure}
\centering
\includegraphics[width=0.75\textwidth]{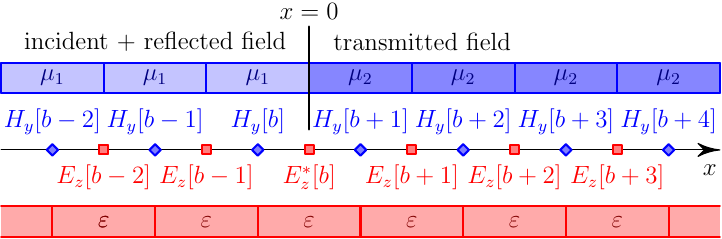}
\caption{FDTD-model of a planar interface between two magnetics}
\label{fig:Magnetics-Interface}
\end{figure}

Now, according to the scheme, in general analogous to the previous subparagraph, we study the interface of two magnetics with relative permeabilities~$\mu_1$ and~$\mu_2$.
Also, for generality, let us assume that the permittivity of both media is equal to~$\varepsilon$, i.e., the contact of pure magnetics for which \mbox{$\varepsilon=1$} is just a special case in our consideration.
This leads us to the following agreement (see Fig.~\ref{fig:Magnetics-Interface})
\begin{equation*}
  \varepsilon_\mathrm{r}[m] \equiv \varepsilon = \const,\quad
  \mu_\mathrm{r}[m] = \begin{cases}
  \mu_1, m\in[0,b],\\
  \mu_2, m\in[b+1,M).
  \end{cases}
\end{equation*}

\begin{remark}
Here, as in the previous subparagraph, there are no additional requirements on the values of~$\mu_1$, $\mu_2$ and~$\varepsilon$ with respect to the Condition~\ref{cond:Perfect-Media}.
So, while the case of \mbox{$\mu_1,\mu_2>\varepsilon$} is primarily meant, \mbox{$\mu_1,\mu_2<\varepsilon$} is not excluded, nor are any intermediate situations.
This, of course, may be a bit contrary to common terminology, but is nevertheless consistent with the models we have adopted, shown in Fig.~\ref{fig:Dielectrics-Interface} and~\ref{fig:Magnetics-Interface}.
\end{remark}

The analysis of the model shown in Fig.~\ref{fig:Magnetics-Interface} indicates that in this case the interface between Medium~1 and Medium~2 is exactly at the node~$E_z^*[b]$.
The condition of continuity of the tangential component of the electric field written for this node (the first equation in~\eqref{eq:Boundary-Conditions}) leads to the first FDTD boundary condition for magnetics
\begin{equation}
\label{eq:Continity-Condition-3}
  1 + \widetilde{r} = \widetilde{t}.
\end{equation}
Similarly~\eqref{eq:Continity-Condition-1}, this expression was obtained after substituting the fields~\eqref{eq:FDTD-Incident-Fields}\,--\,\eqref{eq:FDTD-Transmitted-Fields} into the~\eqref{eq:Boundary-Conditions} and canceling~$E_0e^{i\omega q\Delta_t}$.

For the same reasons as in the previous subparagraph, there is no exact discrete analog of the continuity condition of the tangent component of the magnetic field for the magnetics interface (within the framework of the model shown in Fig.~\ref{fig:Magnetics-Interface}).
This issue is solved by the next theorem.

\begin{theorem}
\label{thrm:2}
The FDTD-analog of the boundary condition relating the tangential components of the magnetic field of a plane monochromatic wave at the interface of two magnetics within the model presented in Fig.~\ref{fig:Magnetics-Interface}, is the following equation
\begin{equation}
\label{eq:Continity-Condition-4}
  e^{i\widetilde{k}_1\Delta_x/2} -
  \widetilde{r}e^{-i\widetilde{k}_1\Delta_x/2} =
  \eta_1 \! \left(\frac{e^{-i\widetilde{k}_2\Delta_x/2}}{\eta_2} +
  i\Omega\varepsilon_0\varepsilon\Delta_x\right) \! \widetilde{t}.
\end{equation}
\end{theorem}

\begin{proof}
The general scheme of the proof is similar to the proof of Theorem~\ref{thrm:1}.
However, since it differs in details, we will give it in full.
First of all, write down the FDTD-analog of Amp\`{e}re's law in the form~\eqref{eq:Ampere-Law-Yee} for the boundary node~$E_z^*[b]$
\begin{equation}
\label{eq:Theorem2-Proof-Step-1}
  \varepsilon_0\varepsilon
  \widetilde{\partial}_t E_z^q[b] =
  \widetilde{\partial}_x H_y^q[b].
\end{equation}

Like for the previous theorem, on this step we use the definition~\eqref{eq:D-Operator-Def} of the operator~$\widetilde{\partial}_x$ and the result~\eqref{eq:Fin-Diff-Operator-Dt} of Lemma~\ref{lemma:D-Actions}
\begin{equation}
\label{eq:Theorem2-Proof-Step-2}
  i \Omega \varepsilon_0 \varepsilon E_z^q[b] =
  \frac1{\Delta_x}
  \left(\widehat{S}_x^{\frac12} - \widehat{S}_x^{-\frac12} \right) H_y^q[b].
\end{equation}

Action of the shift operators~$\widehat{S}^{\pm\frac12}_x$ on the~$H_y^q[b]$ in the right-hand side of the~\eqref{eq:Theorem2-Proof-Step-2} is
\begin{equation}
\label{eq:Theorem2-Proof-Step-3}
  \widehat{S}_x^{\frac12} H_y^q[b] =
  e^{-\widetilde{k}_2\Delta_x/2} H_y^{(t)\,q}[b] =
  -e^{-\widetilde{k}_2\Delta_x/2} \, \widetilde{t} \,
  \frac{E_0}{\eta_2} \, e^{i\omega q\Delta_t},
\end{equation}
\begin{equation}
\label{eq:Theorem2-Proof-Step-4}
  \widehat{S}_x^{-\frac12} H_y^q[b] =
  e^{\widetilde{k}_1\Delta_x/2} H_y^{(i)\,q}[b] +
  e^{-\widetilde{k}_1\Delta_x/2} H_y^{(r)\,q}[b] =
  \left(-e^{\widetilde{k}_1\Delta_x/2} + \widetilde{r} \,
  e^{-\widetilde{k}_1\Delta_x/2}\right)
  \frac{E_0}{\eta_1} \, e^{i\omega q\Delta_t}.
\end{equation}
Here, we take into account that~$\widehat{S}_x^{\frac12}H_y^q[b]$ is magnetic part of the transmitted wave~\eqref{eq:FDTD-Transmitted-Fields}, and~$\widehat{S}_x^{-\frac12} H_y^q[b]$ is the sum of the magnetic fields of the incident~\eqref{eq:FDTD-Incident-Fields} and reflected~\eqref{eq:FDTD-Reflected-Fields} waves.
In addition, one uses~\eqref{eq:Sx-St-Plane-Wave} to compute the application of~$\widehat{S}^{\pm\frac12}_x$ operators on these waves.

As a consequence of the boundary condition~\eqref{eq:Continity-Condition-3}, the left-hand side of the equality~\eqref{eq:Theorem2-Proof-Step-2} can be interpreted as either the sum of the incident and reflected waves or the transmitted wave.
Choosing the second of these options and taking into account expression for the electric field of transmitted wave~\eqref{eq:FDTD-Transmitted-Fields}, we write the left part of~\eqref{eq:Theorem2-Proof-Step-2} in the following form:
\begin{equation}
\label{eq:Theorem2-Proof-Step-5}
  i \Omega \varepsilon_0 \varepsilon E_z^q[b] =
  i \Omega \varepsilon_0 \varepsilon \, \widetilde{t} \,
  E_0 \, e^{i\omega q\Delta_t}.
\end{equation}

Substituting now~\eqref{eq:Theorem2-Proof-Step-3}\,--\,\eqref{eq:Theorem2-Proof-Step-5} into~\eqref{eq:Theorem2-Proof-Step-2}, reducing the common multiplier~$E_0e^{i\omega q\Delta_t}$, multiplying by~$\eta_1\Delta_x$ and grouping the terms with coefficient~$\widetilde{t}$ in the right-hand side of the equality, we finally obtain the~\eqref{eq:Continity-Condition-4}.
\end{proof}

\section{FDTD analog of the Fresnel coefficients}
\label{sect:FDTD-Fresnel-Coeffs}

First, we use Theorems~\ref{thrm:1} and~\ref{thrm:2}, formulated and proved in the previous section, to obtain the FDTD analogs of the classic Fresnel coefficients~\eqref{eq:t-r-Exact}.

\subsection{Dielectric/dielectric interface}

Let us start with the generalized dielectric interface, whose FDTD model is shown in Fig.~\ref{fig:Dielectrics-Interface}.

\begin{theorem}\label{thrm:3}
The FDTD analogs of the Fresnel coefficients~\eqref{eq:t-r-Exact} for the dielectric interface (see Fig.~\ref{fig:Dielectrics-Interface}) are
\begin{equation}
\label{eq:FDTD-t-r-dielectrics}
  \widetilde{t} =
  \frac{2\eta_2\cos\left(\frac{\widetilde{k}_1\Delta_x}{2}\right)}
  {\eta_2\cos\left(\frac{\widetilde{k}_2\Delta_x}{2}\right) +
  \eta_1\cos\left(\frac{\widetilde{k}_1\Delta_x}{2}\right)},\quad
  \widetilde{r} =
  \frac{\eta_2\cos\left(\frac{\widetilde{k}_2\Delta_x}{2}\right) -
  \eta_1\cos\left(\frac{\widetilde{k}_1\Delta_x}{2}\right)}
  {\eta_2\cos\left(\frac{\widetilde{k}_2\Delta_x}{2}\right) +
  \eta_1\cos\left(\frac{\widetilde{k}_1\Delta_x}{2}\right)}.
\end{equation}
\end{theorem}

\begin{proof}
To prove~\eqref{eq:FDTD-t-r-dielectrics}, solve the system of equations~\eqref{eq:Continity-Condition-1} and~\eqref{eq:Continity-Condition-2} for~$\widetilde{r}$ and~$\widetilde{t}$.
To do this, multiply the equality~\eqref{eq:Continity-Condition-1} by $-\eta_1 e^{-i\widetilde{k}_1\Delta_x/2}$, resulting in
\begin{equation}
\label{eq:Thrm3-Step1}
  e^{-i\widetilde{k}_1\Delta_x/2} -
  \widetilde{r} e^{-i\widetilde{k}_1\Delta_x/2} =
  \frac{\eta_1}{\eta_2} e^{-i\widetilde{k}_1\Delta_x/2} \, \widetilde{t}.
\end{equation}

Adding the left- and right-hand sides of~\eqref{eq:Thrm3-Step1} and~\eqref{eq:Continity-Condition-2}, we get an equality that does not contain the coefficient~$\widehat{r}$. Multiplying this expression by~$\eta_2$, we obtain
\begin{equation}
\label{eq:Thrm3-Step2}
  \eta_2 \left( e^{i\widetilde{k}_1\Delta_x/2} +
  e^{-i\widetilde{k}_1\Delta_x/2} \right) =
  \left( \eta_1 e^{-i\widetilde{k}_1\Delta_x/2} +
  \eta_2 e^{-i\widetilde{k}_2\Delta_x/2} +
  i \Omega\mu_0\mu\Delta_x \right) \widetilde{t}.
\end{equation}

Now transform the third term in parentheses on the right-hand side of~\eqref{eq:Thrm3-Step2} using the definition~\eqref{eq:Fin-Diff-Operator-Dt} from Lemma~\ref{lemma:D-Actions} and the dispersion relation~\eqref{eq:Yee-Dispersion-Eq} of Lemma~\ref{lemma:FDTD-Dispersion} written for Medium~2
\begin{equation}
\label{eq:Thrm3-Step3}
  i \Omega\mu_0\mu\Delta_x =
  i \frac{2}{\Delta_t} \frac{S_\mathrm{c}}{n_{\mathrm{r}2}}
  \sin \left(\frac{\widetilde{k}_2\Delta_x}{2} \right) \mu_0\mu\Delta_x.
\end{equation}

Then, using the definitions of Courant number~\eqref{eq:Courant-Number} and relative refractive index~\eqref{eq:Refraction-Index}, the relation~\eqref{eq:Light-Speed} of light speed with vacuum permittivity and permeability, and the definition of wave impedance~\eqref{eq:Impedance}, the coefficients in~\eqref{eq:Thrm3-Step3} can be written
\begin{equation}
\label{eq:Thrm3-Step4}
  \frac{S_\mathrm{c}}{n_{\mathrm{r}2}}
  \frac{\Delta_x}{\Delta_t} \mu_0 \mu =
  \frac{c}{\sqrt{\varepsilon_2\mu}} \mu_0 \mu =
  \sqrt{\frac{\mu_0}{\varepsilon_0}}
  \sqrt{\frac{\mu}{\varepsilon_2}} = \eta_2.
\end{equation}

The sine on the right-hand side of~\eqref{eq:Thrm3-Step3} can be expressed using Euler's formula as
\begin{equation}
\label{eq:Thrm3-Step5}
  \sin \left(\frac{\widetilde{k}_2\Delta_x}{2} \right) =
  \frac{e^{i\widetilde{k}_2\Delta_x/2}-e^{-i\widetilde{k}_2\Delta_x/2}}{2i}.
\end{equation}

Substituting~\eqref{eq:Thrm3-Step4}--\eqref{eq:Thrm3-Step5} into~\eqref{eq:Thrm3-Step3} yields~$\widetilde{t}$ from~\eqref{eq:Thrm3-Step2} as
\begin{equation}
\label{eq:Thrm3-Step6}
  \widetilde{t} = \frac{2\eta_2\cos\left(\frac{\widetilde{k}_1\Delta_x}{2}\right)}
  {\eta_1 e^{-i\widetilde{k}_1\Delta_x/2} + \eta_2 e^{i\widetilde{k}_2\Delta_x/2}}.
\end{equation}
Here, Euler's formula has also been used in writing the numerator.

The real part of the denominator of~\eqref{eq:Thrm3-Step6}
\begin{equation*}
  \re \left( \eta_1 e^{-i\widetilde{k}_1\Delta_x/2} +
  \eta_2 e^{i\widetilde{k}_2\Delta_x/2} \right) =
  \eta_1 \cos\left(\frac{\widetilde{k}_1\Delta_x}{2}\right) +
  \eta_2 \cos\left(\frac{\widetilde{k}_2\Delta_x}{2}\right)
\end{equation*}
coincides with the denominator of the fraction, defining~$\widetilde{t}$ in~\eqref{eq:FDTD-t-r-dielectrics}.
Therefore, it remains only to prove that the imaginary part of this denominator is zero.
This follows from the dispersion relation~\eqref{eq:Yee-Dispersion-Eq} in Lemma~\ref{lemma:FDTD-Dispersion}
\begin{equation}
\label{eq:Thrm3-Step7}
\begin{aligned}
  \im \left( \eta_1 e^{-i\widetilde{k}_1\Delta_x/2} +
  \eta_2 e^{i\widetilde{k}_2\Delta_x/2} \right) &=
  -\eta_1 \sin\left(\frac{\widetilde{k}_1\Delta_x}{2}\right) +
  \eta_2 \sin\left(\frac{\widetilde{k}_2\Delta_x}{2}\right) =\\
  &= \left(-\eta_1 \frac{n_{\mathrm{r}1}}{S_\mathrm{c}} +
  \eta_2 \frac{n_{\mathrm{r}2}}{S_\mathrm{c}} \right)
  \sin\left(\frac{\omega\Delta_t}{2}\right).
\end{aligned}
\end{equation}

With~\eqref{eq:Impedance} and~\eqref{eq:Refraction-Index}, it is easy to verify that the expression in brackets~\eqref{eq:Thrm3-Step7}
\begin{equation*}
  -\eta_1 \frac{n_{\mathrm{r}1}}{S_\mathrm{c}} + \eta_2 \frac{n_{\mathrm{r}2}}{S_\mathrm{c}} =
  \frac{\eta_0}{S_\mathrm{c}} \left(-\sqrt{\frac{\mu}{\varepsilon_1}} \sqrt{\varepsilon_1\mu} +
  \sqrt{\frac{\mu}{\varepsilon_2}} \sqrt{\varepsilon_2\mu} \right) = 0
\end{equation*}
is identically zero.
This finally confirms that the~\eqref{eq:Thrm3-Step6} is the first equality~\eqref{eq:FDTD-t-r-dielectrics}.
To complete the proof, it remains to verify that substituting~$\widetilde{t}$ from~\eqref{eq:FDTD-t-r-dielectrics} into~\eqref{eq:Continity-Condition-1} does indeed lead to the second equality~\eqref{eq:FDTD-t-r-dielectrics}.
\end{proof}

\subsection{Magnetic/magnetic interface}

Now consider the generalized magnetic interface, whose FDTD model is shown in Fig.~\ref{fig:Magnetics-Interface}.

\begin{theorem}\label{thrm:4}
The FDTD analogs of the Fresnel coefficients~\eqref{eq:t-r-Exact} for the magnetic interface (see Fig.~\ref{fig:Magnetics-Interface}) are
\begin{equation}
\label{eq:FDTD-t-r-magnetics}
  \widetilde{t} =
  \frac{2\eta_2\cos\left(\frac{\widetilde{k}_1\Delta_x}{2}\right)}
  {\eta_2\cos\left(\frac{\widetilde{k}_1\Delta_x}{2}\right) +
  \eta_1\cos\left(\frac{\widetilde{k}_2\Delta_x}{2}\right)},\quad
  \widetilde{r} =
  \frac{\eta_2\cos\left(\frac{\widetilde{k}_1\Delta_x}{2}\right) -
  \eta_1\cos\left(\frac{\widetilde{k}_2\Delta_x}{2}\right)}
  {\eta_2\cos\left(\frac{\widetilde{k}_1\Delta_x}{2}\right) +
  \eta_1\cos\left(\frac{\widetilde{k}_2\Delta_x}{2}\right)}.
\end{equation}
\end{theorem}

\begin{proof}
To prove~\eqref{eq:FDTD-t-r-magnetics}, solve~\eqref{eq:Continity-Condition-3}--\eqref{eq:Continity-Condition-4} for~$\widetilde{r}$ and~$\widetilde{t}$.
To do this, multiply the equality~\eqref{eq:Continity-Condition-3} by~$e^{-i\widetilde{k}_1\Delta_x/2}$ and add it with~\eqref{eq:Continity-Condition-4}.
The result of these steps is the equation
\begin{equation}
\label{eq:Thrm4-Step1}
  2 \cos \frac{\widetilde{k}_1\Delta_x}2 = \left(
  \frac{\eta_1}{\eta_2} e^{-i\widetilde{k}_2\Delta_x/2} +
  e^{-i\widetilde{k}_1\Delta_x/2} +
  i\eta_1\Omega\varepsilon_0\varepsilon\Delta_x \right) \widetilde{t}.
\end{equation}
Here, Euler's formula has also been used on the left-hand side.

Now we transform the third summand in parentheses in the right-hand side~\eqref{eq:Thrm4-Step1} using the definition~\eqref{eq:Fin-Diff-Operator-Dt} of~$\Omega$, the FDTD variant of the dispersion relation~\eqref{eq:Yee-Dispersion-Eq}, the definition of impedance~\eqref{eq:Impedance} and the expression for the relative refractive index~\eqref{eq:Refraction-Index}, written for Medium~1
\begin{equation}
\label{eq:Thrm4-Step2}
\begin{aligned}
  i\eta_1\Omega\varepsilon_0\varepsilon\Delta_x &=
  i\eta_1 \frac2{\Delta_t} \sin\left(\frac{\omega\Delta_t}2\right)
  \varepsilon_0\varepsilon\Delta_x =
  i\eta_1 \frac2{\Delta_t} \frac{S_\mathrm{c}}{n_{\mathrm{r1}}}
  \sin\left(\frac{\widetilde{k}_1\Delta_x}2\right)
  \varepsilon_0\varepsilon\Delta_x =\\
  &= i2 \underbrace{\sqrt{\frac{\mu_0\mu_1}{\varepsilon_0\varepsilon}}
  \frac{\Delta_x}{\Delta_t} \frac{S_\mathrm{c}}{\sqrt{\mu_1\varepsilon}}
  \varepsilon_0\varepsilon}_{=1}
  \sin\left(\frac{\widetilde{k}_1\Delta_x}2\right) =
  e^{i\widetilde{k}_1\Delta_x/2} - e^{-i\widetilde{k}_1\Delta_x/2}.
\end{aligned}
\end{equation}
In addition, when writing~\eqref{eq:Thrm4-Step2}, the expression for the speed of light in vacuum~\eqref{eq:Light-Speed}, the definition of Courant number~\eqref{eq:Courant-Number}, and Euler's formula for the sine transformation in the last step are implicitly used.

Substituting the result of transformations~\eqref{eq:Thrm4-Step2} into~\eqref{eq:Thrm4-Step1} gives the following fraction for the coefficient~$\widetilde{t}$
\begin{equation}
\label{eq:Thrm4-Step3}
  \widetilde{t} =
  \frac{2\eta_2\cos\left(\frac{\widetilde{k}_1\Delta_x}{2}\right)}
  {\eta_2 e^{i\widetilde{k}_1\Delta_x/2} +
  \eta_1 e^{-i\widetilde{k}_2\Delta_x/2}}.
\end{equation}

Let us now consider the real and imaginary parts of the denominator of the given fraction
\begin{equation}
\label{eq:Thrm4-Step4}
  \re \left( \eta_2 e^{i\widetilde{k}_1\Delta_x/2} +
  \eta_1 e^{-i\widetilde{k}_2\Delta_x/2} \right) =
  \eta_2\cos\left(\frac{\widetilde{k}_1\Delta_x}{2}\right) +
  \eta_1\cos\left(\frac{\widetilde{k}_2\Delta_x}{2}\right),
\end{equation}
\begin{equation}
\label{eq:Thrm4-Step5}
\begin{aligned}
  \im \left( \eta_2 e^{i\widetilde{k}_1\Delta_x/2} +
  \eta_1 e^{-i\widetilde{k}_2\Delta_x/2} \right) &=
  \eta_2\sin\left(\frac{\widetilde{k}_1\Delta_x}{2}\right) -
  \eta_1\sin\left(\frac{\widetilde{k}_2\Delta_x}{2}\right) =\\
  &= \underbrace{\left(\eta_2 \frac{n_{\mathrm{r}1}}{S_\mathrm{c}} -
  \eta_1 \frac{n_{\mathrm{r}2}}{S_\mathrm{c}}\right)}_{=0}
  \sin\left(\frac{\omega\Delta_t}2\right) = 0.
\end{aligned}
\end{equation}
Here, the dispersion relation~\eqref{eq:Yee-Dispersion-Eq} was also used to write the last equality.
The proof that the expression in brackets in~\eqref{eq:Thrm4-Step5} is zero is analogous to that in Theorem~\ref{thrm:3}.
\begin{equation*}
  \eta_2 \frac{n_{\mathrm{r}1}}{S_\mathrm{c}} -
  \eta_1 \frac{n_{\mathrm{r}2}}{S_\mathrm{c}} =
  \frac{\eta_0}{S_\mathrm{c}} \left(
  \sqrt{\frac{\mu_2}{\varepsilon}}\sqrt{\mu_1\varepsilon} -
  \sqrt{\frac{\mu_1}{\varepsilon}}\sqrt{\mu_2\varepsilon}\right) = 0.
\end{equation*}

Thus, the expressions~\eqref{eq:Thrm4-Step4} and~\eqref{eq:Thrm4-Step5} prove that the equality~\eqref{eq:Thrm4-Step3} is the first equation in~\eqref{eq:FDTD-t-r-magnetics}.
Its substitution into the boundary condition~\eqref{eq:Continity-Condition-3} allows us to express the coefficient~$\widetilde{r}$ in the form of the second equation in~\eqref{eq:FDTD-t-r-magnetics}, which completes the proof.
\end{proof}

\begin{remark}
\label{note:FDTD-Fresnel-Coeffs}
Note that $\widetilde{k}_{(1,2)}\Delta_x/2$ in~\eqref{eq:FDTD-t-r-dielectrics} and~\eqref{eq:FDTD-t-r-magnetics} must be expressed in terms of the incident wave frequency using~\eqref{eq:Yee-Dispersion-Eq} from Lemma~\ref{lemma:FDTD-Dispersion}.
Thus, the FDTD variants of the Fresnel coefficients~\eqref{eq:FDTD-t-r-dielectrics} and~\eqref{eq:FDTD-t-r-magnetics} unlike their exact continuous analogs~\eqref{eq:t-r-Exact} are determined not only by the impedances~$\eta$ of the bordering media, but also by the incident wave frequency, as well as by specific parameters of simulation, namely---the ratio between the Courant number~$S_\mathrm{c}$ and the relative refractive indices~$n_\mathrm{r}$ of both media.
\end{remark}

\begin{proposition}\label{prop:FDTD-Exact-Compliance}
The results~\eqref{eq:FDTD-t-r-dielectrics} and~\eqref{eq:FDTD-t-r-magnetics} of both Theorems we proved earlier in this paragraph converge to the exact expressions~\eqref{eq:t-r-Exact} at the continuous limit transition~$\Delta_x\rightarrow0$.
\end{proposition}

\begin{proof}
This is easily established by acting similarly to Remark~\ref{note:Dispersion}, given after Lemma~\ref{lemma:FDTD-Dispersion}.
\end{proof}

\section{Discussion}
\label{sect:Discussion}

In this section, we discuss the main applications of Theorems~\ref{thrm:3} and~\ref{thrm:4}, and the implications they have for the practice of FDTD modeling.

\subsection{Peculiarities of the FDTD Fresnel coefficients}
\label{sect:FDTD-Fresnel-Coefficients}

As the next step in our investigation, let us examine in detail the properties of the coefficients~$\widetilde{r}$ and~$\widetilde{t}$ that follow from the Theorem~\ref{thrm:3} and~\ref{thrm:4}.

\begin{figure}
\centering
\includegraphics[width=0.95\textwidth]{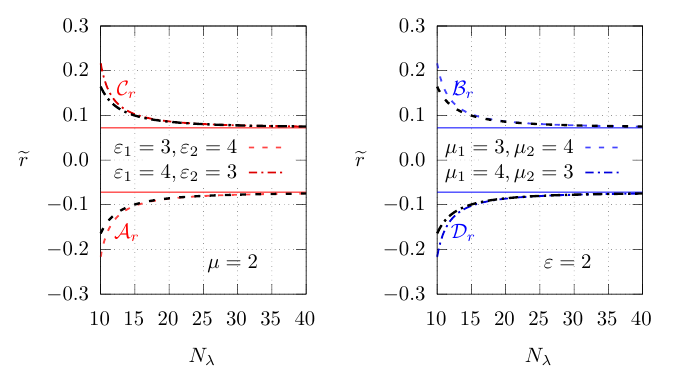}
\includegraphics[width=0.95\textwidth]{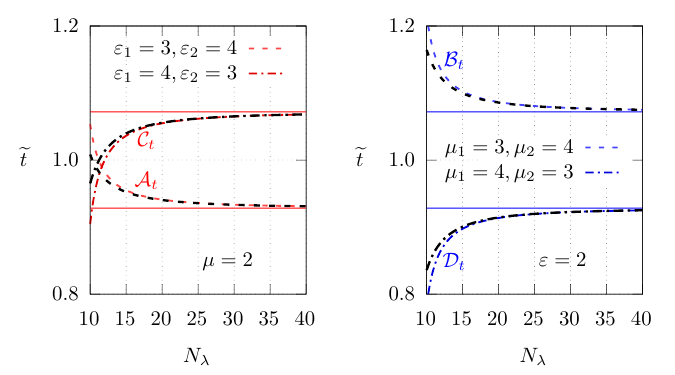}
\caption{FDTD Fresnel coefficients for reflection~$\widetilde{r}$ (two images in the upper row) and transmission~$\widetilde{t}$ (lower row) in the case of the interface between two dielectrics (left column) and magnetics (right column) at slightly different impedances of the two media \mbox{$\eta_1/\eta_2=\sqrt{4/3}\approx1.16$} and \mbox{$\eta_1/\eta_2=\sqrt{3/4}\approx0.87$} (the specific values~$\varepsilon_\mathrm{r}$ and~$\mu_\mathrm{r}$ at which the curves are plotted are shown in the corresponding parts of the figure)}
\label{fig:FDTD-Fresnel-Coeffs}
\end{figure}

So, the characteristic features of the FDTD Fresnel coefficients~\eqref{eq:FDTD-t-r-dielectrics} and~\eqref{eq:FDTD-t-r-magnetics} are clearly demonstrated by Fig.~\ref{fig:FDTD-Fresnel-Coeffs}, constructed for the case of a not very contrasting interface (this means that the impedances of the two media are chosen to differ little from each other).
Thus, for \mbox{$\eta_1/\eta_2\approx1.16$} and \mbox{$\eta_1/\eta_2\approx0.87$} (see details in the caption to Fig.~\ref{fig:FDTD-Fresnel-Coeffs} and in parts of the figure itself), the exact formulas~\eqref{eq:t-r-Exact} predict pairs of values \mbox{$r\approx-0.07$}, \mbox{$t\approx0.93$}, and \mbox{$r\approx0.07$}, \mbox{$t\approx1.07$} regardless of the incident wave frequency.
These values are marked in Fig.~\ref{fig:FDTD-Fresnel-Coeffs} with solid horizontal lines.

At the same time, as already noted in Remark~\ref{note:FDTD-Fresnel-Coeffs}, the results of FDTD calculations~\eqref{eq:FDTD-t-r-dielectrics} and~\eqref{eq:FDTD-t-r-magnetics} depend significantly on the frequency~$\omega$ of the incident wave.
This dependence, as the main feature, is primarily illustrated in Fig.~\ref{fig:FDTD-Fresnel-Coeffs}.
It should be clarified that instead of the frequency~$\omega$ in Fig.~\ref{fig:FDTD-Fresnel-Coeffs} there is a discretization parameter~$N_\lambda$, defined as the number (not necessarily integer) of the Yee grid nodes per wavelength~$\lambda$ of the incident wave in vacuum
\begin{equation}\label{eq:N-lambda}
  \lambda = N_\lambda \Delta_x.
\end{equation}

With~$N_\lambda$, the frequency of the incident wave is expressed according to
\begin{equation}\label{eq:FDTD-Omega-N-lambda}
  \frac{\omega\Delta_t}2 =
  \frac{\pi}{N_\lambda} S_\mathrm{c}.
\end{equation}
As was noted in Remark~\ref{note:FDTD-Fresnel-Coeffs}, the sine of this particular expression according to the formula~\eqref{eq:Yee-Dispersion-Eq} of the Lemma~\ref{lemma:FDTD-Dispersion} with accuracy up to the multiplier~$n_\mathrm{r(1,2)}/S_\mathrm{c}$ determines the values~$\widetilde{k}_{(1,2)}\Delta_x/2$ appearing in the cosines of the formulas~\eqref{eq:FDTD-t-r-dielectrics} and~\eqref{eq:FDTD-t-r-magnetics}.

To obtain~\eqref{eq:FDTD-Omega-N-lambda}, one uses the exact equation~\eqref{eq:Dispersion-Law} of the dispersion of electromagnetic waves in vacuum (\mbox{$n_\mathrm{r}=1$}), the relation of wave number and wavelength \mbox{$k=2\pi/\lambda$} (see~\cite{LandauV8}, \cite{Jackson} or any other electrodynamics textbooks), the expression~\eqref{eq:N-lambda} and the definition of Courant number~\eqref{eq:Courant-Number}
\begin{equation*}
  \frac{\omega\Delta_t}2 =
  \frac{kc\Delta_t}2 =
  \frac{2\pi}{\lambda} \frac{c\Delta_t}{2} =
  \frac{\pi}{N_\lambda} \frac{c\Delta_t}{\Delta_x} =
  \frac{\pi}{N_\lambda} S_\mathrm{c}.
\end{equation*}

With increasing~$N_\lambda$ (which, by virtue of~\eqref{eq:FDTD-Omega-N-lambda}, is equivalent to a decrease in the frequency~$\omega$ of the incident wave), the spatial step~$\Delta_x$ of the Yee grid becomes negligibly small compared to the wavelength~$\lambda$.
This is in full agreement with Proposition~\ref{prop:FDTD-Exact-Compliance} explains the dependence of the coefficients~$\widetilde{r}$ and~$\widetilde{t}$ on~$N_\lambda$, which asymptotically tend to their exact counterparts~$r$ and~$t$ at \mbox{ $N_\lambda\rightarrow\infty$}.
The curves shown in Fig.~\ref{fig:FDTD-Fresnel-Coeffs} clearly demonstrate this feature and, in addition, allow us to make numerical estimates of the errors in determining the Fresnel coefficients by the FDTD method for a particular interfaces of  dielectrics or magnetics at any chosen level of discretization~$N_\lambda$.

Moreover, Fig.~\ref{fig:FDTD-Fresnel-Coeffs} contains information about the two types interfaces of dielectrics and magnetics: the dashed lines represent the dependences on~$N_\lambda$ of the FDTD coefficients~$\widetilde{r}$ and~$\widetilde{t}$ when \mbox{$\varepsilon_1<\varepsilon_2$} (red curves for dielectrics interface) and \mbox{$\mu_1<\mu_2$} (blue lines in the magnetics example), while the point-dashed lines correspond to the opposite situations.

Despite the seemingly diverse variations shown in Fig.~\ref{fig:FDTD-Fresnel-Coeffs} in terms of different values of~$\varepsilon_{(1,2)}$ and~$\mu_{(1,2)}$, many of these results can be generalized using the concept of wave impedance~\eqref{eq:Impedance}.
Thus, it is obvious that the case of the interface of two dielectrics with \mbox{$\varepsilon_1=3$} and \mbox{$\varepsilon_2=4$} (and the same value of \mbox{$\mu=2$}) is equivalent from the point of view of asymptotic values \mbox{$\lim\limits_{N_\lambda\rightarrow\infty}\widetilde{r}=r$} and \mbox{$\lim\limits_{N_\lambda\rightarrow\infty}\widetilde{t}=t$} to the situation of the interface of two magnetics with comparable parameters: \mbox{$\mu_1=4$} and \mbox{$\mu_2=3$} (and the common value of \mbox{$\varepsilon=2$}).
Both situations (such dielectrics and magnetics) correspond to the example with \mbox{$\eta_1/\eta_2\approx1.16$} discussed above and correspond to the case of reflection with half-wave loss (\mbox{$r<0$}).
Similar conclusions can be drawn about the dielectric and magnetic interfaces at the same value of \mbox{$\eta_1/\eta_2\approx0.87$}.

\begin{remark}
A remark should be made here on the terminology adopted in optics, when one distinguishes between optically more or less dense media.
It is guided by the fact that in the spectral range of electromagnetic waves corresponding to the visible light, the relative permeability~$\mu_\mathrm{r}$ of the vast majority of media (with the very rare exception of special, artificially created metamaterials) is practically indistinguishable from unity.
In such conditions, in the case of \mbox{$\varepsilon_1<\varepsilon_2$} it is usual to speak about reflection of light from optically denser medium, which is accompanied by half-wave loss.
The half-wave loss is understood as the sign of the coefficient~$r$, which is negative in this situation, and this means that the actual directions of the vectors~$\mathbf{E}^{(r)}$ and~$\mathbf{H}^{(r)}$ are opposite to those depicted in Fig.~\ref{fig:General-Geometry}.
Our study, results of which are the curves in Fig.~\ref{fig:FDTD-Fresnel-Coeffs}, indicates that in the general case, to compute the Fresnel coefficients and the vector directions defined by them~$\mathbf{E}^{(r)}$, $\mathbf{H}^{(r)}$ and~$\mathbf{E}^{(t)}$, $\mathbf{H}^{(t)}$ should be based exactly on the relation between the impedances~$\eta_1/\eta_2$ of the two bordering media.
\end{remark}

Furthermore, it should be noted that in Fig.~\ref{fig:FDTD-Fresnel-Coeffs}, the red dot-dashed line in the left column in the top row is qualitatively identical to the corresponding dashed line in the right column of that row.
The same correspondence holds for the red dashed and blue dot-dashed curves.
Thus, we can conclude that for the same ratios~$\eta_1/\eta_2$, not only do the asymptotic values of the FDTD coefficients~$\widetilde{r}$ coincide, but also the qualitative nature of their dependence on~$N_\lambda$.
However, the same cannot be said about the dependence \mbox{$\widetilde{t}=\widetilde{t}(N_\lambda)$}.
This clearly follows from the nature of the curves in the lower row of Fig.~\ref{fig:FDTD-Fresnel-Coeffs}, which show that for \mbox{$\eta_1/\eta_2>1$}, in the case of dielectrics \mbox{$\widetilde{t}>t$}, while for magnetics with a similar ratio~$\eta_1/\eta_2$, the opposite is true: \mbox{$\widetilde{t}<t$}.
When moving to the case \mbox{$\eta_1/\eta_2<1$}, everything for \mbox{$\widetilde{t}=\widetilde{t}(N_\lambda)$} changes to the exact opposite.
All these observations can be generalized in the form of the following three Corollaries.

\begin{corollary}\label{Cor:1}
The asymptotic values \mbox{$\lim\limits_{N_\lambda\rightarrow\infty}\widetilde{r}=r$} and \mbox{$\lim\limits_{N_\lambda\rightarrow\infty}\widetilde{t}=t$} remain the same for the boundary between both dielectrics and magnetics, and are determined only by the specific value of the ratio~$\eta_1/\eta_2$.
\end{corollary}

\begin{proof}
The proof of this Corollary is based on Proposition~\ref{prop:FDTD-Exact-Compliance} and the details of the models of the boundaries between dielectrics and magnetics presented in Fig.~\ref{fig:Transition-Layers}.
The numerical parameters of permittivities and permeabilities shown in Fig.~\ref{fig:Transition-Layers} are chosen to be the same as those in the corresponding parts of Fig.~\ref{fig:FDTD-Fresnel-Coeffs}.
In this regard, the boundary model shown in Fig.~\ref{fig:Transition-Layers}{\em a)} determines the behavior of the families of curves~$\mathcal{A}_r$ and~$\mathcal{A}_t$ shown in Fig.~\ref{fig:FDTD-Fresnel-Coeffs}.
The same analogy is accepted between the remaining parts of Fig.~\ref{fig:Transition-Layers} and the families of curves~$\mathcal{B}$, $\mathcal{C}$ and~$\mathcal{D}$ shown in Fig.~\ref{fig:FDTD-Fresnel-Coeffs}.

\begin{figure}
\centering
\includegraphics[width=0.85\textwidth]{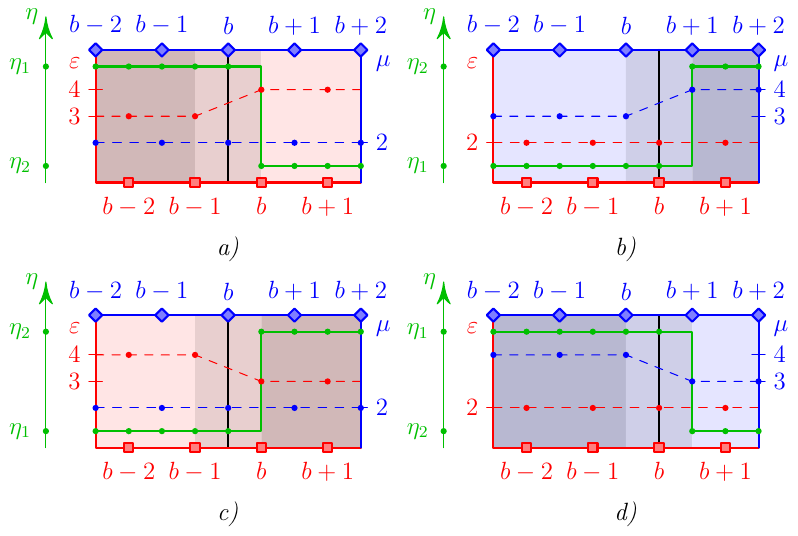}
\caption{Formation features of a transition layer between two dielectrics: {\em a)}~\mbox{$\eta_1>\eta_2$}, {\em c)}~\mbox{$\eta_1<\eta_2$}; and for magnetics: {\em b)}~\mbox{$\eta_1<\eta_2$}, {\em d)}~\mbox{$\eta_1>\eta_2$}}
\label{fig:Transition-Layers}
\end{figure}

The Yee grid nodes located at the exact boundary between the two media in Fig.~\ref{fig:Transition-Layers}, in full accordance with the FDTD models shown in Fig.~\ref{fig:Dielectrics-Interface} and Fig.~\ref{fig:Magnetics-Interface}, are marked with the index~$b$.
However, Fig.~\ref{fig:Transition-Layers} clearly shows the main feature of the FDTD method, which cannot be seen in Fig.~\ref{fig:Dielectrics-Interface} and Fig.~\ref{fig:Magnetics-Interface}.
This feature is that the Yee grid nodes closest to the exact boundary and having numbers~$b-1$, $b$ (for dielectrics) and~$b$, $b+1$ (for magnetics), due to the discontinuous behavior of the material parameters of the boundary media, form a transition layer with a length of~$\Delta_x$.

Throughout this transition layer, the impedance of the Yee grid changes with some ``delay'' in two ways, shown in Fig.~\ref{fig:Transition-Layers} by green lines.
The nature of this change in the transition layer of dielectrics with \mbox{$\varepsilon_1<\varepsilon_2$} is identical to the corresponding change in impedance in the transition layer of magnetics with \mbox{$\mu_1>\mu_2$} (compare the green curves in Fig.~\ref{fig:Transition-Layers}{\em a)} and Fig.~\ref{fig:Transition-Layers}{\em d)}).
A similar picture is obtained when comparing the curves in Fig.~\ref{fig:Transition-Layers}{\em c)} and Fig.~\ref{fig:Transition-Layers}{\em b)}).

To prove this Corollary, it is essential that the range of values \mbox{$N_\lambda\rightarrow\infty$}, as noted above, is equivalent to such a choice of~$\Delta_x$ that can be considered negligible compared to~$\lambda$.
This means that the sizes of the transition regions in Fig.~\ref{fig:Transition-Layers} at \mbox{$N_\lambda\rightarrow\infty$} turn out to be negligible, and the nature of the impedance changes shown in Fig.~\ref{fig:Transition-Layers} for both dielectrics and magnetics becomes practically indistinguishable from that for the point boundary model.
This proves that the asymptotic values of the FDTD Fresnel coefficients are determined only by the ratio~$\eta_1/\eta_2$, regardless of the type of boundary media.
\end{proof}

\begin{corollary}\label{Cor:2}
For the interfaces of both dielectrics and magnetics, the nature of the dependence of the FDTD-coefficients~$\widetilde{r}$ on~$N_\lambda$ is identical at the same ratio~$\eta_1/\eta_2$.
It can be argued that when $\eta_1>\eta_2$, the FDTD variant of the Fresnel coefficient for reflection is always underestimated compared to its exact value \mbox{$\widetilde{r}<r$}, and in the reverse situation \mbox{$\eta_1<\eta_2$}, the opposite is true: \mbox{$\widetilde{r}>r$}.
\end{corollary}

\begin{corollary}\label{Cor:3}
For the dependence of the FDTD coefficients~$\widetilde{t}$ on~$N_\lambda$, there is a lot of variation here.
Thus, when \mbox{$\eta_1>\eta_2$}, for dielectrics there is an excess \mbox{$\widetilde{t}>t$}, while for magnetics reduced values \mbox{$\widetilde{t}<t$} are observed.
When transitioning to the case $\eta_1<\eta_2$, everything changes to the exact opposite.
\end{corollary}

Now let us discuss the peculiarities of FDTD-calculating Fresnel coefficients, which are due to the variability of modeling parameters.
As already noted in the introduction to this article, the Courant number~$S_\mathrm{c}$, defined according to~\eqref{eq:Courant-Number}, plays a crucial role in this context.

Thus, the colored curves (red for dielectrics and blue for magnetics) in Fig.~\ref{fig:FDTD-Fresnel-Coeffs} are plotted for standard conditions, which imply that the Courant number \mbox{$S_\mathrm{c}=1$}.
This mode is optimal only when modeling electromagnetic processes in a vacuum (\mbox{$n_\mathrm{r}=1$}).
This is because in this situation, the speed of electromagnetic waves is exactly equal to~$c$, which means that when moving from the current iteration of the Yee algorithm~\eqref{eq:FDTD-Faraday-Law}\,--\,\eqref{eq:Source-Ez} to the next (after a time interval~$\Delta_t$), changes in the electromagnetic field in space can only propagate effectively if the distance between the nodes of the Yee grid is \mbox{$\Delta_x=c\Delta_t$} (which corresponds to the choice \mbox{$S_\mathrm{c}=1$}).

In the paper~\cite{Makarov2024a}, it was shown that the choice \mbox{$S_\mathrm{c}=1$} ceases to be optimal in the case of FDTD modeling of electromagnetic processes in non-dispersive homogeneous media with \mbox{$n_\mathrm{r}\neq1$}.
The reason for this is the effects associated with the numerical dispersion of the FDTD method, described in Lemma~\ref{lemma:FDTD-Dispersion}.
Moreover, in a situation where \mbox{$n_\mathrm{r}<1$} FDTD modeling performed according to the standard scheme (\mbox{$S_\mathrm{c}=1$}) is in principle impossible, since the implementation of the calculation algorithm~\eqref{eq:FDTD-Faraday-Law}\,--\,\eqref{eq:Source-Ez} very quickly leads to a divergent solution.

According to the work~\cite{Makarov2024a}, the solution to these problems is the choice of~\eqref{eq:Magic-Courant-Number}, which reconciles the propagation speed of harmonic waves with the parameters of the Yee grid.
However, in the context of current paper, this choice needs to be refined, since, due to the very formulation of the problem, we are dealing with a boundary between two media with different values of~$n_\mathrm{r}$ in general case.

\begin{definition}\label{def:FDTD-Optimal-Mode}
Let us call the iterative algorithm Yee~\eqref{eq:FDTD-Faraday-Law}\,--\,\eqref{eq:Source-Ez}, performed for models of the interface between two media, shown in Fig.~\ref{fig:Dielectrics-Interface} or~\ref{fig:Magnetics-Interface}, in which the Courant number is set equal to
\begin{equation}\label{eq:Optimal-Courant-Number}
  S_\mathrm{c} = \min(n_{\mathrm{r}1}, n_{\mathrm{r}2}),
\end{equation}
the optimal mode of FDTD simulation.
\end{definition}

\begin{remark}
This definition can be generalized to the case of a larger number of media \mbox{$N>2$}, but this is not necessary for the purposes of the current study.
\end{remark}

\begin{remark}\label{rem:FDTD-Optimal-Mode}
Choosing the minimum of the two values in~\eqref{eq:Optimal-Courant-Number} guarantees that the implementation of the computational scheme~\eqref{eq:FDTD-Faraday-Law}\,--\,\eqref{eq:Source-Ez} in the Conditions~\ref{cond:Perfect-Media} and~\ref{cond:Currents-and-Charges-Absence} we have adopted will not lead to a divergent solution for dielectric and magnetic cases.
At the same time, it is obvious that the optimal mode does not completely eliminate the effects caused by the numerical dispersion of FDTD.
\end{remark}

The black curves in Fig.~\ref{fig:FDTD-Fresnel-Coeffs} are plotted for the situation when FDTD modeling is performed in the optimal mode~\eqref{eq:Optimal-Courant-Number}.
It can be seen that the positive effect on the accuracy of the calculations of the coefficients~$\widetilde{r}$ and~$\widetilde{t}$ is clearly evident in the region of very coarse discretization (when~$N_\lambda$ is sufficiently small), i.e., for relatively high-frequency waves.
In the opposite case of slowly changing electromagnetic fields, the difference in the calculation of Fresnel coefficients~\eqref{eq:FDTD-t-r-dielectrics} and~\eqref{eq:FDTD-t-r-magnetics} in standard and optimal modes becomes insignificant.
We will return to a more detailed discussion of this issue later in connection with the analysis of Fig.~\ref{fig:R-T-Coeffs-Errors-S}.

\subsection{FDTD version of the reflection and transmission coefficients}
\label{sect:FDTD-Reflection-and-Transmission-Coefficients}

Now, after a detailed examination of the FDTD Fresnel coefficients, let us move on to the analysis of the FDTD analogues of the reflection and transmission coefficients.
These quantities are easy to introduce based on the definition of their exact analogues~\eqref{eq:R-T-Exact}:
\begin{equation}
  \label{eq:FTDT-R-T}
  \widetilde{R} = \widetilde{r}\,^2,\quad
  \widetilde{T} = \frac{\eta_1}{\eta_2} \widetilde{t}\,^2.
\end{equation}

\subsubsection{Low impedance contrast interface}

\begin{figure}
\centering
\includegraphics[width=0.5\textwidth]{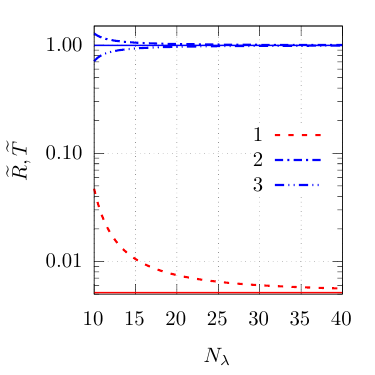}
\caption{FDTD-coefficients for reflection~$\widetilde{R}$ (curves~1) and transmission~$\widetilde{T}$ (curves~2 and~3) in the case of the interface between two media with slightly different impedances: \mbox{$\eta_1/\eta_2=\sqrt{4/3} \approx 1.16$} (curves~2 for dielectrics and~3---for magnetics) and \mbox{$\eta_1/\eta_2=\sqrt{3/4} \approx 0.87$} (curves~3 for dielectrics and~2---for magnetics)}
\label{fig:R-T-Coeffs}
\end{figure}

Fig.~\ref{fig:R-T-Coeffs} shows the dependence of the coefficients~\eqref{eq:FTDT-R-T} on the discretization parameter~$N_\lambda$ for the interface between two dielectrics and two magnetics in the case of a small difference between their impedances.
The same parameter values were used to construct this figure as were used in Fig.~\ref{fig:FDTD-Fresnel-Coeffs}.
We will return to the case of a high-contrast interface later when discussing Fig.~\ref{fig:R-T-Coeffs-High-Contrast}.

First, note that Fig.~\ref{fig:R-T-Coeffs} confirms the well-known~\cite{LandauV8,Jackson} fact that at the boundary between any media with weak impedance contrast, the reflection is relatively small (\mbox{$R\approx0$}), since the overwhelming majority of the electromagnetic wave energy passes practically unhindered from the first medium to the second (\mbox{$T\approx1$}).

Secondly, Fig.~\ref{fig:R-T-Coeffs} shows that the FDTD method in the case of coarse discretization (with relatively small~$N_\lambda$, i.e. for high-frequency waves), leads to a very significant difference between the reflection coefficients~$\widetilde{R}$ and transmission coefficients~$\widetilde{T}$ and the corresponding exact values, observed during numerical simulation.
The qualitative nature of these differences can be formulated in the form of two corollaries.

\begin{corollary}\label{Cor:4}
The values of the FDTD reflection coefficients~$\widetilde{R}$ are always overestimated compared to the exact values~$R$, regardless of the specific simulation conditions.
\end{corollary}

\begin{proof}
The validity of this fact follows from the direct application of the first formula~\eqref{eq:FTDT-R-T} to all curves shown in the two images in the upper row of Fig.~\ref{fig:FDTD-Fresnel-Coeffs}, the nature of which was discussed in Corollary~\ref{Cor:2}.
\end{proof}

\begin{corollary}\label{Cor:5}
The values of FDTD transmission coefficients~$\widetilde{T}$ are overestimated compared to their exact counterparts~$T$ in the case of a boundary between dielectrics with \mbox{$\eta_1>\eta_2$} or magnetics with \mbox{$\eta_1<\eta_2$}.
In the opposite cases, the FDTD-transmission are underestimated \mbox{$\widetilde{T}<T$}.
\end{corollary}

\begin{proof}
The proof of this result is based on applying the second formula~\eqref{eq:FTDT-R-T} to all curves shown in both images in the bottom row of Fig.~\ref{fig:FDTD-Fresnel-Coeffs}.
All data required for this argument can be conveniently presented in the form of Table~\ref{tab:FDTD-Transmission-Errors}.
\begin{table}[t]
\caption{Data from Fig.~\ref{fig:FDTD-Fresnel-Coeffs}, necessary for proving Corollary~\ref{Cor:5}}
\label{tab:FDTD-Transmission-Errors}
\center
\begin{tabular}{@{}cccc@{}}
\toprule
Curves & Type of Interface & Features of~$t$ and~$\widetilde{t}$ & Results for~$\widetilde{T}$\\
\midrule
$\mathcal{A}_t$ & dielectrics with \mbox{$\eta_1>\eta_2$} & \mbox{$t<1$}, \mbox{$\widetilde{t}>t$} & \multirow{2}*{\mbox{$\widetilde{T}>T$}}\\
$\mathcal{B}_t$ & magnetics with \mbox{$\eta_1<\eta_2$} & \mbox{$t>1$}, \mbox{$\widetilde{t}>t$} & \\
$\mathcal{C}_t$ & dielectrics with \mbox{$\eta_1<\eta_2$} & \mbox{$t>1$}, \mbox{$\widetilde{t}<t$} & \multirow{2}*{\mbox{$\widetilde{T}<T$}}\\
$\mathcal{D}_t$ & magnetics with \mbox{$\eta_1>\eta_2$} & \mbox{$t<1$}, \mbox{$\widetilde{t}<t$} & \\
\bottomrule
\end{tabular}
\end{table}

The third column of this table contains information from Corollary~\ref{Cor:3} about the relationship between the values~$\widetilde{t}$ and~$t$, obtained from the analysis of Fig.~\ref{fig:FDTD-Fresnel-Coeffs}.
The fourth column shows the result of applying formula~\eqref{eq:FTDT-R-T} to the data in the third column.
Combining the contents of the second and fourth columns clearly represents the main result of this Corollary.
\end{proof}

Thirdly, since the ordinate axis in Fig.~\ref{fig:R-T-Coeffs} is plotted on a logarithmic scale, this allows us not only to see the dependence of the absolute values of the coefficients~$\widetilde{R}$ and~$\widetilde{T}$ on~$N_\lambda$, but also to qualitatively estimate their relative errors.
Thus, it is clear from Fig.~\ref{fig:R-T-Coeffs} that the deviation of curves~$1$, $2$, and~$3$ from the corresponding asymptotes is described by the following inequality: $\left|\widetilde{T}-T\right|<\left|\widetilde{R}-R\right|$.
This means that the accuracy of calculating the intensity of the transmitted wave using the FDTD method in this case of weak contrast between the impedances of the two media is higher than that of the reflected wave.

To quantitatively describe the accuracy of the FDTD method at the interface, it is convenient to use the following values, which characterize the relative errors in the estimation of the reflection and transmission coefficients~\eqref{eq:FTDT-R-T}:
\begin{equation}
  \label{eq:FTDT-delta-R-T}
  \delta_R = \frac{\left|\widetilde{R} - R\right|}{R} \cdot 100\%,\quad
  \delta_T = \frac{\left|\widetilde{T} - T\right|}{T} \cdot 100\%.
\end{equation}

\begin{figure}
\centering
\includegraphics[width=0.85\textwidth]{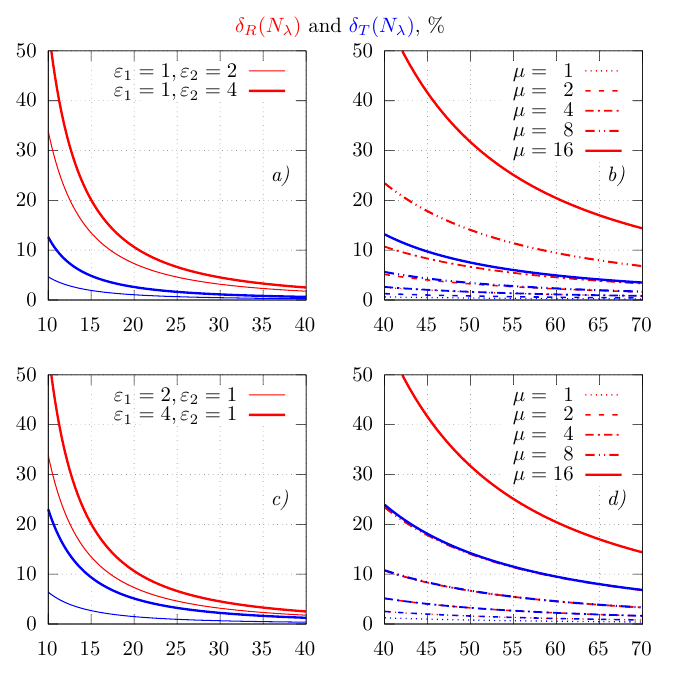}
\caption{Influence on the errors of FDTD reflection coefficients (red curves) and transmission coefficients (blue lines) of dielectric permittivities (left column of images) and magnetic permeabilities (right column) in the case of a boundary between two dielectrics with \mbox{$\eta_1>\eta_2$} (top row of images) and \mbox{$\eta_1<\eta_2$} (bottom row).
When constructing images~{\em a)} and~{\em c)}, the total magnetic permeability was chosen to be unity \mbox{$\mu=1$}; for case {\em b)} the values \mbox{$\varepsilon_1=1$} and \mbox{$\varepsilon_2=4$} are given; in case {\em d)} \mbox{$\varepsilon_1=4$} and \mbox{$\varepsilon_2=1$}.
Thus, the thick solid lines shown on the left side of the figure continuously transition into thin dotted curves on the right side, covering two ranges of the discretization parameter: \mbox{$N_\lambda\in[10,40]$} and \mbox{$N_\lambda\in[40,70]$}.
The Courant number \mbox{$S_\textrm{c}=1$} was used in the construction of all parts of this figure.
All other numerical parameters used in the construction of the image are indicated on the corresponding parts of the figure}
\label{fig:R-T-Coeffs-Errors}
\end{figure}

The nature of the dependence of the relative errors~\eqref{eq:FTDT-delta-R-T} of the FDTD reflection and transmission coefficients on the discretization parameter~$N_\lambda$ in the case of the interface between two dielectrics is shown in Fig.~\ref{fig:R-T-Coeffs-Errors}.
This figure clearly demonstrates the fact that the dependencies~$\delta_R(N_\lambda)$ and~$\delta_T(N_\lambda)$ are monotonically decreasing functions.
The numerical values of these functions are significantly influenced by all material parameters (permittivities and permeability) of both boundary media.
It should be emphasized here that the influence of the parameter~$\mu$ on the nature of the curves~$\delta_R(N_\lambda)$ and~$\delta_T(N_\lambda)$ in the dielectric interface model shown in Fig.~\ref{fig:Dielectrics-Interface} is not ``intuitively obvious'', but it turns out to be very significant, which is well reflected in the right column of images in Fig.~\ref{fig:R-T-Coeffs-Errors}.

Moving on to the analysis of Fig.~\ref{fig:R-T-Coeffs-Errors}, we first note that in the case of weak contrast between the impedances of the boundary media (\mbox{$\eta_1/\eta_2\approx1$}), the relative errors of the FDTD reflection coefficients always exceed the corresponding errors of the transmission coefficients: \mbox{$\delta_R>\delta_T$}.
Here, however, it should be noted that this is not always true when switching to the opposite situation: \mbox{$\eta_1/\eta_2\gg1$}, as well as \mbox{$\eta_1/\eta_2\ll1$} (see below for analysis of Fig.~\ref{fig:R-T-Coeffs-High-Contrast}).

Secondly, an increase in the contrast of the impedances of the boundary media corresponds to an increase in both~$\delta_R$ and~$\delta_T$ (see the transition from thin to thick lines in the images in the left column of Fig.~\ref{fig:R-T-Coeffs-Errors}).
Another characteristic feature of the~$\delta_R$ and~$\delta_T$, which can be seen on the left side of Fig.~\ref{fig:R-T-Coeffs-Errors}, is that the error~$\delta_R$ in determining the FDTD reflection coefficient is invariant with respect to the substitution \mbox{$\eta_1\leftrightarrow\eta_2$}\label{page:ImpedanceSymmetry}, while~$\delta_T$ is determined not only by the values of the impedances~$\eta_1$ and~$\eta_2$, but also by the order in which they occur.
Thus, when a wave passes through the boundary between two dielectrics with \mbox{$\eta_1>\eta_2$}, the error~$\delta_T$ in determining the FDTD transmission coefficient is always smaller than in the opposite case \mbox{$\eta_1<\eta_2$}.
The latter fact follows from a comparison of the blue curves in Fig.~\ref{fig:R-T-Coeffs-Errors}{\em a)} and~\ref{fig:R-T-Coeffs-Errors}{\em c)}.

Thirdly, the accuracy of determining FDTD reflection and transmission coefficients depends not only on the ratio of the impedances of two dielectrics, $\eta_1/\eta_2$, but also on the value of their total magnetic permeability~$\mu$.
This fact follows from the analysis of the curves shown in the right column of the images in Fig.~\ref{fig:R-T-Coeffs-Errors}.
Thus, since in the model of the interface between two dielectrics shown in Fig.~\ref{fig:Dielectrics-Interface}, the permeability of both media is the same (\mbox{$\mu=\const$}), it is obvious that the impedance ratio~$\eta_1/\eta_2$ does not depend on it.
This means that for all types of lines shown in Fig.~\ref{fig:R-T-Coeffs-Errors}{\em b)} the ratio~$\eta_1/\eta_2$ is exactly the same as for the curves shown by thick lines in Fig.~\ref{fig:R-T-Coeffs-Errors}{\em a)}.
The same applies to Fig.~\ref{fig:R-T-Coeffs-Errors}{\em d)} and~\ref{fig:R-T-Coeffs-Errors}{\em c)}.
Analysis of the nature of the curves~$\delta_R(N_\lambda)$ and~$\delta_T(N_\lambda)$ curves shown on the right side of Fig.~\ref{fig:R-T-Coeffs-Errors} allows us to conclude that an increase in~$\mu$, all other conditions being equal, increases the errors in determining both the reflection coefficient and the transmission coefficient.
In this case, the increase in~$\mu$ in the given case of weak contrast between the impedances of the two media increases~$\delta_R$ more than~$\delta_T$.

\begin{remark}
The case of the boundary between magnetics with permeability values equivalent to those shown in Fig.~\ref{fig:R-T-Coeffs-Errors} does not require separate detailed consideration, since the corresponding curves~$\delta_R(N_\lambda)$ and~$\delta_T(N_\lambda)$ are practically identical to those already shown.
The only difference is that the upper row of images in Fig.~\ref{fig:R-T-Coeffs-Errors} for equivalent magnets will correspond to the situation \mbox{$\eta_1<\eta_2$}, while the lower row will correspond to the case \mbox{$\eta_1>\eta_2$}.
\end{remark}

\begin{figure}
\centering
\includegraphics[width=0.85\textwidth]{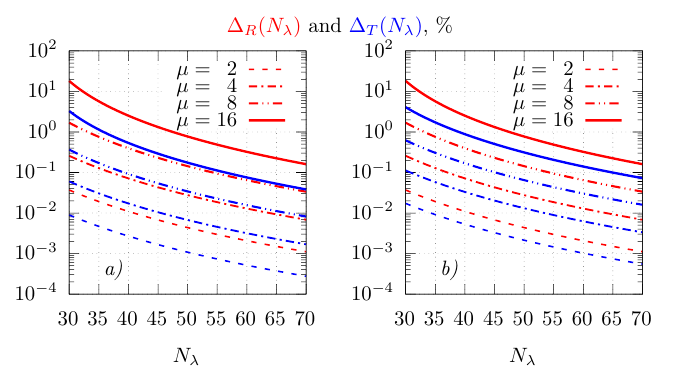}
\caption{The difference between the errors of the FDTD coefficients for reflection~$\Delta_R$ (red curves) and transmission~$\Delta_T$ (blue lines), calculated using the standard scheme (\mbox{$S_\mathrm{c}=1$}) and in the optimal mode (\mbox{$S_\mathrm{c}=\min(n_{\mathrm{r}1},n_{\mathrm{r}2})$}) in the case of a boundary between two dielectrics: {\em a)}~\mbox{$\eta_1>\eta_2$} (\mbox{$\varepsilon_1=1$} and \mbox{$\varepsilon_2=4$}); {\em b)}~\mbox{$\eta_1<\eta_2$} (\mbox{$\varepsilon_1=4$} and \mbox{$\varepsilon_2=1$}).
The case with a common permeability \mbox{$\mu=1$} is not shown, since there is no difference between the standard and optimal modeling schemes in this case.}
\label{fig:R-T-Coeffs-Errors-S}
\end{figure}

Now let us turn to the question of the influence of the choice of the Courant number~$S_\mathrm{c}$ on the errors in calculating the reflection~$\delta_R(N_\lambda)$ and the transmission~$\delta_T(N_\lambda)$ coefficient using the FDTD method.
To do this, we will compare the errors in the FDTD calculation of the corresponding coefficients performed using the standard scheme (\mbox{$S_\mathrm{c}=1$}) and in the optimal mode (see Definition~\ref{def:FDTD-Optimal-Mode} at the end of \S\ref{sect:FDTD-Fresnel-Coefficients}).
We denote the differences between the first and second results as~$\Delta_{(R,T)}$, and we present in Fig.~\ref{fig:R-T-Coeffs-Errors-S} the corresponding graphs for the case of dielectrics whose numerical values of material parameters are similar to the more contrasting dielectrics presented in Fig.~\ref{fig:R-T-Coeffs-Errors}.

Fig.~\ref{fig:R-T-Coeffs-Errors-S} shows that setting a reasonable value for the sampling parameter~$N_\lambda$ has a much greater impact on the accuracy of FDTD calculations of reflection and transmission coefficients than choosing the optimal value for the Courant number.
Thus, at \mbox{$N_\lambda=40$}, the maximum benefit from using the optimal simulation scheme with the parameters shown in Fig.~\ref{fig:R-T-Coeffs-Errors-S} is only \mbox{$\Delta_R\approx3\%$} (see the curve plotted at \mbox{$\mu=16$}), while the error in determining the reflection coefficient exceeds \mbox{$\delta_R>50\%$} (see the corresponding curves in Fig.~\ref{fig:R-T-Coeffs-Errors}{\em b)} and~\ref{fig:R-T-Coeffs-Errors}{\em d)}).
Obviously, against the background of such a significant error, a 3\% improvement in calculation accuracy seems insignificant.

In conclusion of this consideration, we also note that the choice of the optimal Courant number can make a significant difference in the calculations (\mbox{$\Delta\gtrsim10\%$}) only in the case of very coarse discretization (see Fig.~\ref{fig:R-T-Coeffs-Errors-S} for the range of values \mbox{$N_\lambda<40$}), i.e. where this difference does not give any real advantage with a fundamentally poor simulation result.
At the same time, as the discretization of the wavelength of the incident wave improves with an increase in~$N_\lambda$, the difference in accuracy between the standard and optimal FDTD calculations decreases monotonically and very rapidly, amounting to \mbox{$\Delta_R<0.2\%$} (in the case of \mbox{$\mu=16$}) already at \mbox{$N_\lambda=70$} (which is still poor, since the error~$\delta_R$ in this case is approximately~15\%).

\subsubsection{Interface with high contrast of impedances}

\begin{figure}
\centering
\includegraphics[width=0.85\textwidth]{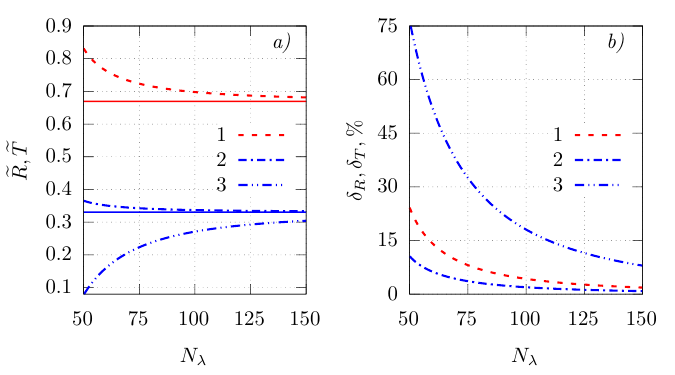}
\caption{{\em a)} FDTD coefficients for reflection~$\widetilde{R}$ (curves~1) and transmission~$\widetilde{T}$ (curves~2 and~3) in the case of the interface between two media with a high difference in impedances: \mbox{$\eta_1/\eta_2=10$} (curves~2 for dielectrics and~3 for magnetics) and \mbox{$\eta_1/\eta_2~=~0.1$} (curves~3 for dielectrics and~2 for magnetics); {\em b)} Relative errors of the corresponding coefficients under the same simulating conditions}
\label{fig:R-T-Coeffs-High-Contrast}
\end{figure}

Now let us consider the case of a more contrasting interface.
Fig.~\ref{fig:R-T-Coeffs-High-Contrast}{\em a)} shows the dependence of coefficients~\eqref{eq:FTDT-R-T} on the sampling parameter~$N_\lambda$ for two media with impedance ratios \mbox{$\eta_1/\eta_2=10$} and \mbox{$\eta_1/\eta_2=0.1$}: dielectrics with \mbox{$\varepsilon_{(1,2)}=\{1,100\}$} and \mbox{$\mu=2$}, as well as their equivalent magnetics.
This figure should be analyzed in conjunction with Fig.~\ref{fig:R-T-Coeffs}, however, it should be noted that the scale on the vertical axis in Fig.~\ref{fig:R-T-Coeffs} is logarithmic, while in Fig.~\ref{fig:R-T-Coeffs-High-Contrast}{\em a)}---it is linear.
This is due to the fact that the deviations of all values of coefficients~$\widetilde{R}$ and~$\widetilde{T}$ from their exact analogues in Fig.~\ref{fig:R-T-Coeffs-High-Contrast}{\em a)} are within the same order of magnitude, while in Fig.~\ref{fig:R-T-Coeffs} there are different scales of these changes.

First, we note that in this case, reflection prevails over transmission for both the exact values of the coefficients \mbox{$R>T$} and their FDTD analogues \mbox{$\widetilde{R}>\widetilde{T}$}.
This, of course, was to be expected for the case of a boundary between materials with a large difference in wave impedances, but it is worth emphasizing that in order to observe this, FDTD modeling should be performed in the region of significantly large values of the discretization parameter~$N_\lambda$.
Thus, in Fig.~\ref{fig:R-T-Coeffs}, the entire range of values~$N_\lambda$ is limited to the interval \mbox{$[10,40]$}, while for the medium impedance ratios shown in Fig.~\ref{fig:R-T-Coeffs-High-Contrast}, this range is in principle inaccessible for FDTD modeling.
The latter is due to the fact that the significant difference between the medium impedance~$\eta$ and the vacuum impedance~$\eta_0$ is due to large values of the relative refractive index of the medium~$n_\mathrm{r}$, which entails a reduction in the wavelength in the medium (see formulas~\eqref{eq:Impedance}, \eqref{eq:Refraction-Index}, and~\eqref{eq:N-lambda}).
Accordingly, to demonstrate the features of FDTD modeling of reflection and transmission coefficients at a high-contrast interface, the range~$N_\lambda$ in Fig.~\ref{fig:R-T-Coeffs-High-Contrast}{\em a)} was shifted to the region of large values \mbox{$[50,150]$}.

Secondly, Fig.~\ref{fig:R-T-Coeffs-High-Contrast} clearly demonstrates the misconception of the naive assumption that could have arisen earlier when analyzing Fig.~\ref{fig:R-T-Coeffs}\,--\,\ref{fig:R-T-Coeffs-Errors-S}, that the relative error in determining the reflection coefficient in the case of a well-reflecting interface (\mbox{$R>T$}) should be less than the relative error in the transmission coefficient \mbox{$\delta_R<\delta_T$}.
In reality, however, in Fig.~\ref{fig:R-T-Coeffs-High-Contrast}{\em a)} the deviation of curve~$1$ from its asymptote is smaller than the corresponding deviation only when compared to curve~$3$, but not~$2$.

Thus, before comparing the errors~$\delta_R$ and~$\delta_T$ with each other, we should first consider the asymmetry of curves~$2$ and~$3$ relative to their asymptote.
In fact, we have already dealt with this effect earlier (see the blue curves in Fig.~\ref{fig:R-T-Coeffs-Errors}{\em a)} and~\ref{fig:R-T-Coeffs-Errors}{\em c)}, as well as Fig.~\ref{fig:R-T-Coeffs-Errors-S}{\em a)} and~\ref{fig:R-T-Coeffs-Errors-S}{\em b)} and the discussion regarding the replacement of \mbox{$\eta_1\leftrightarrow\eta_2$} on page~\pageref{page:ImpedanceSymmetry}), but then it was not so evident.
Here, existence of this asymmetry is much more clear, since, along with the impedances, the relative refractive indices of the two media also differ significantly.
This explains why the accuracy of the FDTD calculation is the worst for curve~$3$ in Fig.~\ref{fig:R-T-Coeffs-High-Contrast}{\em a)}.
This is due to the fact that in the case of curve~$3$, the first medium in which the incident wave propagates from the source to the interface is so optically dense that its amplitude and shape by the time it reaches the interface are already very distorted compared to the original ones due to the numerical dispersion of the FDTD method with insufficient discretization~$N_\lambda$.
Of course, in this case, there is no point in talking about any accurate determination of the amplitude and intensity of the wave that has passed into the second medium.
Note that choosing the optimal Courant number~\eqref{eq:Optimal-Courant-Number} cannot correct this problem, since the minimum refractive index in this case is determined by the value~$n_{\mathrm{r}2}$ for the second medium.

Now, after explaining the nature of such strong asymmetry between curves~$2$ and~$3$ in Fig.~\ref{fig:R-T-Coeffs-High-Contrast}{\em a)}, we can compare the errors~$\delta_R$ and~$\delta_T$ in the case of a boundary between media with high impedance contrast.
A perfect illustration of this issue is given in Fig.~\ref{fig:R-T-Coeffs-High-Contrast}{\em b)}, which should be compared with Fig.~\ref{fig:R-T-Coeffs-Errors}.
The main conclusion of this comparison can be formulated as follows: the relative errors in determining the FDTD transmission coefficient, when a reasonable modeling scheme is chosen (i.e., when \mbox{$n_{\mathrm{r}1}<n_{\mathrm{r}2}$}), are always smaller than the corresponding errors for the reflection coefficient \mbox{$\delta_T<\delta_R$}.
The only exception to this rule is the situation of a sufficiently high-contrast interface between two media when \mbox{$n_{\mathrm{r}1}>n_{\mathrm{r}2}$}.

And in conclusion of the analysis of the case of a boundary between media with high impedance contrast, we note that in this situation the difference between the standard and optimal simulation modes is negligible.
This follows from the fact that the maximum value of the parameter~$\Delta$ for all curves shown in Fig.~\ref{fig:R-T-Coeffs-High-Contrast}{\em b)} does not exceed the value \mbox{$\Delta_{T3}<0.1\%$}.

\subsection{Possible topics for further research}
\label{sect:Further-Research}

Earlier in this paragraph, we discussed such features of FDTD modeling of normal incidence of harmonic plane electromagnetic waves on a plane interface of linear homogeneous and isotropic dielectrics and magnetics as the behavior of Fresnel coefficients~$\widetilde{r}$ and~$\widetilde{t}$, the reflection~$\widetilde{R}$ and the transmission~$\widetilde{T}$ coefficients, their asymptotic values, as well as the relative errors~$\delta_R$ and~$\delta_T$ depending on the characteristics of the boundary media (permittivity, permeability and impedance), as well as such modeling parameters as the wavelength discretization parameter~$N_\lambda$ and the Courant number~$S_\mathrm{c}$.

Thus, this paper describes in detail a large number of results obtained at the current stage of work.
At the same time, research in this area can be continued further, since there are many interesting options for its development.
Below, we will limit ourselves to listing only a few possible promising aims that can be set in subsequent works.
\begin{enumerate}
  \item Direct analytical proof is required for Corollary~\ref{Cor:2} and~\ref{Cor:3}, the validity of which was established in this work solely on the basis of observations and analysis of the results of numerical calculations presented in Fig.~\ref{fig:FDTD-Fresnel-Coeffs}.
Analytical proofs of qualitative error behavior for interface-resolving schemes are increasingly valued in the numerical analysis community; cf. \cite{Wang2025} for a recent example in the context of high-order methods for wave equations.
  \item There is independent interest to consider cases with \mbox{$\varepsilon_\mathrm{r}<1$} and/or \mbox{$\mu_\mathrm{r}<1$}, which do not fall outside the scope of Condition~\ref{cond:Perfect-Media}, but have not been specifically considered in this work.
This may be useful when simulating electromagnetic processes occurring in such plasma~\cite{Pavlenko2019}, for which the layered medium model is valid.
  \item From the point of view of practical applications in the field of information transmission in the form of complex electromagnetic signals, the analysis of the application of the results of this study beyond the harmonic approximation is promising.
  \item One of the most important applied problems requiring separate detailed consideration and going beyond the scope of the approximation~\ref{cond:Perfect-Media} is obtaining analogues of formulas~\eqref{eq:FDTD-t-r-dielectrics} and~\eqref{eq:FDTD-t-r-magnetics} taking into account energy dissipation (i.e., in the case of \mbox{$\varepsilon_\mathrm{r},\mu_\mathrm{r}\in\mathbb{C}$}).
  \item Most of the content in \S\ref{sect:FDTD-Boundary-Conditions}\,--\,\S\ref{sect:Discussion} is based on the medium boundary models shown in Figs.~\ref{fig:Dielectrics-Interface} and~\ref{fig:Magnetics-Interface}.
In this regard, an important development of this work is the generalization of Theorems~\ref{thrm:1}\,--\,\ref{thrm:4} to the case of media differing simultaneously in values~$\varepsilon_\mathrm{r}$ and~$\mu_\mathrm{r}$.
Note that this generalization is possible both within the scope of Condition~\ref{cond:Perfect-Media} and outside its limitations.
  \item Another example of problems related to research beyond Condition~\ref{cond:Perfect-Media} is the study of the features of FDTD modeling of reflection and transmission of electromagnetic waves at interfaces when at least one of the media is left-handed media~\cite{Veselago1967,Pendry2004}.
\end{enumerate}

We would like to highlight point 5 in the list above.
Although real-world interfaces typically involve simultaneous changes to both~$\varepsilon_\mathrm{r}$ and $\mu_\mathrm{r}$, we deliberately chose not to address this case in this paper, as it requires more careful analysis.

\begin{figure}
\centering
\includegraphics[width=0.85\textwidth]{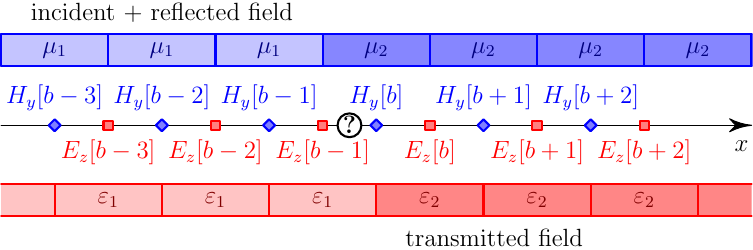}
\caption{FDTD-model of ``real interfaces'' involve simultaneous changes in both~$\varepsilon_\mathrm{r}$ and~$\mu_\mathrm{r}$}
\label{fig:Real-Interface-FDTD-Model}
\end{figure}

This is due to two reasons.
Firstly, the question arises of a rigorous proof of the equivalence of choosing~$E_z[b-1]$ or~$H_y[b]$ of the Yee grid as the exact boundary between the nodes, as shown in Fig.~\ref{fig:Real-Interface-FDTD-Model}.
This is not a particularly complex task, but its straightforward solution requires discussing a vast number of technical details, which would greatly increase the length of this article and become an unnecessary distraction.
However, the following rather obvious result can be anticipated.

\begin{hypothesis}
The asymptotic behavior of the Fresnel coefficients, as well as the reflection and transmission coefficients in the \mbox{$N_\lambda\rightarrow\infty$} region for the FDTD interface model shown in Fig.~\ref{fig:Real-Interface-FDTD-Model} also tends to be their exact analogs to the same extent as for the models shown in Fig.~\ref{fig:Dielectrics-Interface} and~\ref{fig:Magnetics-Interface}.
Moreover, numerical error remains comparable and Corollaries~\ref{Cor:1}, \ref{Cor:2} and~\ref{Cor:4} qualitatively persist when impedance contrast dominates.
\end{hypothesis}

Secondly, taking into account a more general interface model will require clarification of the Definition~\ref{def:FDTD-Optimal-Mode}, since even the equality of the optical refractive indices of the adjacent layers of the dielectric and magnetic material does not in all cases guarantee the receipt of a stable solution (contrary to the Remark~\ref{rem:FDTD-Optimal-Mode}).
This can be easily understood by considering, for example, the interfaces between materials with \mbox{$\varepsilon_{\mathrm{r}1}=2$}, \mbox{$\mu_{\mathrm{r}1}=1$} \big/ \mbox{$\varepsilon_{\mathrm{r}2}=1$}, \mbox{$\mu_{\mathrm{r}2}=2$} (and the opposite), using the transition-layer idea demonstrated in Fig.~\ref{fig:Transition-Layers}.
Investigating the details of this issue takes us away from the main objectives of this article, which is devoted to investigating the features of the FDTD Fresnel coefficients.

That is is why point 5 is a promising direction for further research.

\section{Conclusion}

This article thoroughly examines the peculiarities of FDTD simulations of normal incidence of harmonic plane electromagnetic waves on a planar interface between linear, homogeneous, isotropic, and perfect dielectrics and magnetics, explicitly accounting for their permittivity and permeability.

At the beginning of the study, specific conditions were set, within which all subsequent arguments were carried out.
Section~\ref{sect:Electrodynamics-Basics} establishes exact analytical expressions for Fresnel coefficients~\eqref{eq:t-r-Exact} and reflection/transmission coefficients~\eqref{eq:R-T-Exact}, serving as benchmarks for FDTD accuracy evaluation.

Sections~\ref{sect:FDTD-Basics}--\ref{sect:FDTD-Dispersion} outline the FDTD methodology, including Yee update equations~\eqref{eq:FDTD-Faraday-Law}--\eqref{eq:FDTD-Ampere-Law}, TFSF source implementation~\eqref{eq:Source-Hy}--\eqref{eq:Source-Ez}, and the dispersion relation~\eqref{eq:Yee-Dispersion-Eq}.

The main research part of the work is concentrated further in~\S\ref{sect:FDTD-Boundary-Conditions}, \S\ref{sect:FDTD-Fresnel-Coeffs} and~\S\ref{sect:Discussion}.
In~\S\ref{sect:FDTD-Boundary-Conditions} two basic FDTD models are defined for the interfaces between dielectrics and magnetics, within the scope of which two theorems on the corresponding boundary conditions are formulated and proven.
Subsection \S\ref{sect:FDTD-DefinitionsAndLemmas} contains necessary definitions and lemmas, which are used further to obtain the main results of this section.
In~\S\ref{sect:DielectricsBoundary}, it is proven that~\eqref{eq:Continity-Condition-1} and~\eqref{eq:Continity-Condition-2} are FDTD analogues of the boundary conditions relating the tangential components of the electromagnetic field of a plane monochromatic wave at the interface of two dielectrics.
The validity of~\eqref{eq:Continity-Condition-3} and~\eqref{eq:Continity-Condition-4}, which are FDTD analogues of the boundary conditions at the interface of two magnetics, is proven in~\S\ref{sect:MagneticsBoundary}.

In \S\ref{sect:FDTD-Fresnel-Coeffs} two more theorems about FDTD analogs of the Fresnel coefficients at the interface between dielectrics and magnetics are formulated and proven.

\S\ref{sect:Discussion} contains a detailed discussion of the results of applying Theorems~\ref{thrm:3} and~\ref{thrm:4} in the practice of FDTD simulation.
In particular, it is established that the results~\eqref{eq:FDTD-t-r-dielectrics} and~\eqref{eq:FDTD-t-r-magnetics} of both Theorems~\ref{thrm:3} and~\ref{thrm:4} remain the same for the boundary between dielectrics and magnetics, and are determined only by the specific value of the ratio~$\eta_1/\eta_2$ and converge to the exact expressions~\eqref{eq:t-r-Exact} in the continuous limit as \mbox{$\Delta_x\rightarrow0$}.
When discussing this result, a transition layer model was formulated, shown in Fig.~\ref{fig:Transition-Layers}, which shows that the impedance of the Yee grid changes at the interface between any media with delay.
Furthermore, it was established that for the interfaces of both dielectrics and magnetics, the nature of the dependence of the FDTD coefficients~$\widetilde{r}$ on~$N_\lambda$ is identical at the same ratio~$\eta_1/\eta_2$, while for the dependence of the FDTD coefficients~$\widetilde{t}$ on~$N_\lambda$, there is a lot of variation.
Next, at the end of~\S\ref{sect:FDTD-Fresnel-Coefficients}, the definition of the optimal mode of FDTD simulation is formulated, based on a specific choice of the Courant number.

After that, at the beginning of~\S\ref{sect:FDTD-Reflection-and-Transmission-Coefficients}, FDTD analogues of the reflection and transmission coefficients are introduced, the behavior and errors of which are discussed further in two opposite cases: for weak and high-contrast (in terms of impedance differences) interfaces between dielectrics and magnetics.
Corollaries~\ref{Cor:4}--\ref{Cor:5} establish that FDTD reflection coefficients are always overestimated (\mbox{$\widetilde{R}>R$}), while transmission coefficients show direction-dependent bias.
FDTD errors~$\delta_R$ and~$\delta_T$ increase not only with rising contrast of the impedances of the boundary media, but also with increasing values of their common permeability (in the case of dielectrics) or permeability (in the case of magnetics).
It is also shown that the value of~$\delta_T$ is affected by the order of the impedances of the boundary media, while the error~$\delta_R$ is invariant with respect to the replacement \mbox{$\eta_1\leftrightarrow\eta_2$}.
Furthermore, in~\S\ref{sect:FDTD-Reflection-and-Transmission-Coefficients} it is established that the relative errors~\eqref{eq:FTDT-delta-R-T} of the FDTD reflection coefficients always exceed the corresponding errors of the transmission coefficients \mbox{$\delta_R>\delta_T$} in the case of weak contrast between the impedances of the boundary media.
However, this is not always true when the situation is reversed.
It has also been found that adjusting the sampling parameter~$N_\lambda$ has a much greater effect on the accuracy of FDTD calculations of reflection and transmission coefficients than choosing the optimal value of the Courant number~$S_\mathrm{c}$.
To conclude this paragraph, \S\ref{sect:Further-Research} lists possible directions for further research in this area.

Thus, this work develops and tests an approach to analyzing the accuracy of the FDTD method in computational electrodynamics, which, using discrete grids and approximations, can be used both to reproduce fundamental physical results and to solve problems for which there are currently no analytical solutions.
The results of this study may be useful in the design of photonic, nano-optical, and radio engineering devices for the development of correct simulating techniques and their verification.

Beyond the specific results presented, this work underscores the continued relevance of rigorous accuracy analysis for classical discretizations, even in an era of advanced structure-preserving methods.
By quantifying interface-induced errors in the Yee scheme and interpreting them through a transition-layer model, we provide tools that are immediately applicable to industrial simulation, educational instruction, and the validation of next-generation Maxwell solvers.
The explicit error estimates derived in this work may serve as benchmarks for validating higher-dimensional FDTD implementations and structure-preserving discretizations of Maxwell's equations.

\section*{Acknowledgments}

The authors express their sincere gratitude to John B. Schneider, Associate Professor of Washington State University for useful discussion of FDTD dispersion effects and all anonymous reviewers for their valuable suggestions that improved this work.

The study was performed within the framework of a state assignment of the FRC Komi SC UB RAS (topic No. 125020501562-1), and also with financial support from the Ministry of Science and Higher Education of Russia within Agreement No. 075‑15‑2025‑455 dated 26.05.2025 during the conduct of NEXAFS studies.

\bibliographystyle{elsarticle-num}

\bibliography{references}

\end{document}